\documentclass{article}

\pdfoutput=1

\usepackage{guit}
\usepackage{arxiv}
\usepackage{graphicx}
\usepackage{amsmath}
\usepackage{amssymb}
\usepackage{mathtools}
\usepackage{indentfirst}
\usepackage{caption}
\usepackage{subcaption}
\usepackage{multirow}
\usepackage[utf8]{inputenc} %
\usepackage[T1]{fontenc}    %
\usepackage{hyperref}       %
\usepackage{url}            %
\usepackage{booktabs}       %
\usepackage{amsfonts}       %
\usepackage{nicefrac}       %
\usepackage{microtype}      %
\usepackage{lipsum}
\usepackage[algo2e]{algorithm2e}
\usepackage{algorithm}
\usepackage{algpseudocode}
\usepackage{float}
\usepackage{upgreek,textgreek}

\DeclarePairedDelimiter{\ceil}{\lceil}{\rceil}

\usepackage[backend=bibtex]{biblatex}
\usepackage{textgreek}
\bibliography{references}

\title{Non intrusive reduced order modeling of parametrized PDEs by kernel POD and neural networks}

\author{
  Matteo Salvador \\
  MOX, Department of Mathematics\\
  Politecnico di Milano\\
  P.zza Leonardo da Vinci 32, 20133 Milan, Italy \\
  \texttt{matteo1.salvador@polimi.it} \\
  \And
  Luca Dede' \\
  MOX, Department of Mathematics\\
  Politecnico di Milano\\
  P.zza Leonardo da Vinci 32, 20133 Milan, Italy \\
  \texttt{luca.dede@polimi.it} \\
   \And
  Andrea Manzoni \\
  MOX, Department of Mathematics\\
  Politecnico di Milano\\
  P.zza Leonardo da Vinci 32, 20133 Milan, Italy \\
  \texttt{andrea1.manzoni@polimi.it} \\
}

\begin{document}
\maketitle

\begin{abstract}
We propose a nonlinear reduced basis method for the efficient approximation of parametrized partial differential equations (PDEs), exploiting kernel proper orthogonal decomposition (KPOD) for the generation of a reduced-order space and neural networks for the evaluation of the reduced-order approximation. In particular, we use KPOD in place of the more classical POD, on a set of high-fidelity solutions of the problem at hand to extract a reduced basis. This method provides a more accurate approximation of the snapshots' set featuring a lower dimension, while maintaining the same efficiency as POD. %
A neural network (NN) is then used to find the coefficients of the reduced basis by following a supervised learning paradigm and shown to be effective in learning the map between the time/parameter values and the projection of the high-fidelity snapshots onto the reduced space.
In this NN, both the number of hidden layers and the number of neurons vary according to the intrinsic dimension of the differential problem at hand and the size of the reduced space. This adaptively built NN attains good performances in both the learning and the testing phases.
Our approach is then tested on two benchmark problems, a one-dimensional wave equation and a two-dimensional nonlinear lid-driven cavity problem. %
We finally compare the proposed KPOD-NN technique with a POD-NN strategy, showing that KPOD allows a reduction of the number of modes that must be retained to reach a given accuracy in the reduced basis approximation. For this reason, the NN built to find the coefficients of the KPOD expansion is smaller, easier and less computationally demanding to train than the one used in the POD-NN strategy.
\end{abstract}

\keywords{Reduced order modeling \and Kernel proper orthogonal decomposition \and Proper orthogonal decomposition \and Neural networks \and Parametrized PDEs} %

\section{Introduction}
\label{section: Introduction}
Reduced order modeling (ROM) techniques represent a very efficient approach for the numerical approximation of problems involving the repeated solution of differential equations arising from engineering and applied sciences  \cite{Fresca1, Hesthaven1, Maulik, Pagani, Karniadakis1, Regazzoni, San}. Their aim is to replace the original large-dimension numerical problem, which is typically called either high-fidelity approximation or full order model (FOM), by a reduced problem of substantially smaller dimension, to provide a fast and reliable approximation to the PDE solution for each newly queried parameter instance.

Depending on the context, there are several strategies to generate the reduced problem from the high-fidelity one.
The approach followed by the reduced basis (RB) method consists in the projection of the high-fidelity problem upon a low-dimensional subspace made of specially selected basis functions, built from  a set of high-fidelity solutions corresponding to suitably chosen parameters, e.g., through proper orthogonal decomposition (POD) \cite{Quarteroni2}. The dimension of the reduced model should be in principle fairly lower than the FOM one. POD is an example of linear dimensionality reduction technique. Such a technique might
show poor performances when dealing with highly nonlinear PDEs, because of the need of a huge number of POD modes to reach a desired accuracy of the reduced order approximation. Nonlinear dimensionality reduction techniques, such as kernel proper orthogonal decomposition (KPOD), have become popular also in the field of reduced order modeling for parametrized PDEs due to their better capability of spanning low-dimensional nonlinear trial manifolds \cite{Xing2015, Xing2016}. By exploiting the better representation power of nonlinear maps, these approaches  condense better information from the underlying PDE keeping an even smaller reduced dimension than POD \cite{Scholkopf}.  However, generating the reduced order approximation efficiently when dealing with nonlinear dimensionality reduction techniques might be far from being trivial. KPOD can be either used to perform the forward mapping, i.e. nonlinear dimensionality reduction, or the backward mapping, which is also called {\em pre-image} reconstruction \cite{Diez2020}. Focusing on the forward mapping, KPOD projects the high-fidelity solutions by means of a nonlinear map to a high dimensional space (called {\em feature space}), where we are likely to obtain linear separability: then, by performing a dimensionality reduction, KPOD maps these data onto a reduced space that is in principle very similar to the space spanned by the most relevant modes in POD \cite{Scholkopf}. By employing a suitable {\em kernel trick} -- as usually done when dealing with kernel principal component analysis (KPCA) \cite{Scholkopf} -- we do not need to build the (very high-dimensional) feature space  explicitly, finally relying on a suitably modified version of the  snapshots correlation (or Gram) matrix, involving the evaluation of a bilinear kernel function instead of classical inner products. An important aspect, in this respect, is to find a suitable strategy to map the reduced solution back onto the high-fidelity space -- a task that, in the case of POD, would simply entail the left-multiplication of the reduced approximation by the matrix collecting the vectors representing the POD modes \cite{Quarteroni2}. In the case of KPOD, for instance, this task would involve the minimization of a discrepancy functional. The main motivations rely on the cost entailed by this stage, and the impossibility to provide a Galerkin projection to generate the (low-dimensional) reduced-order problem, which is solved for any new parameter instance.

For all these reasons, we propose to rely on neural network (NN) regression to determine, in a non intrusive way, the reduced order approximation for any new parameter instance, once a reduced basis has been built through KPOD.
Thanks to their high representational power and flexibility \cite{Barron, Cybenko, Kolmogorov, Zhang}, neural networks are increasingly employed in the numerical approximation of PDEs \cite{Fresca1, Hesthaven1, Karniadakis3, Karniadakis4, Regazzoni}.
Several recent works have shown possible applications of neural networks to parametrized PDEs -- thanks to their approximation capabilities, their extremely favorable computational performances during online testing phases, and their relative easiness of implementation --  both from a theoretical \cite{Jentzen, bhattacharya2020model, Petersen2020, Petersen, Shen} and a computational standpoint \cite{Fresca1, Hesthaven1, Mucke2020, Regazzoni}. A priori characterizations of the NN architecture complexity in terms of the accuracy of the reconstructed solution have been reported in \cite{Petersen}, showing that the NN complexity should scale with the intrinsic dimension of the PDE solution manifold, i.e. with the size of the reduced basis extracted from the PDE under investigation to properly capture its behavior. From a computational standpoint, feedforward NNs and autoencoders, have been employed in several ways to determine the reduced order approximation in a data-driven and less intrusive way (avoiding, e.g., the computational cost entailed by projection-based ROMs), but still relying on a linear trial manifold built, e.g., through POD. In \cite{guo2019data, guo2018reduced, Hesthaven1, San, Hesthaven2} the solution of nonlinear, time-dependent ROMs, for any new parameter instance, has been replaced by the evaluation of NN-based regression models, built using a fixed size feedforward neural network. Few attempts have been made in order to model the reduced order manifold where the approximation is sought (avoiding, e.g., the linear superimposition of POD modes) through NNs, see, e.g., \cite{gonzalez2018deep,lee2020model}. In the framework of deep learning, variational autoencoders have been combined with a fixed size feedforward NN and eventually POD in \cite{Fresca3, Fresca1, Fresca2} to build non intrusive ROMs similarly to \cite{Hesthaven1}, capturing however more details by using the intrinsic dimension of a differential problem for the reduced order approximation, that is, the number of parameters (plus one, to take into account the time coordinate) the solution depends on.

In this work, we propose a non-intrusive ROM technique that combines KPOD with an adaptively built NN, i.e. a NN where the number of layers and neurons is automatically adjusted  to the complexity of the parametrized PDE at hand. This complexity is captured by the dimension of the reduced basis extracted from the FOM snapshots through KPOD. This latter technique  is, at the same time, as efficient as the POD, yet  capable of collecting most of the information in a small number of modes of the kernel matrix, thus performing data compression even more efficiently than POD.

The paper is structured as follows: in Section~\ref{section: KPODNNMethod} we describe our ROM method involving a combined use of KPOD and NN, and we apply it to a general parametrized PDEs. In Section \ref{section: MathematicalModeling} we introduce the mathematical models on which our metholodology will be tested, i.e. wave equation and Navier-Stokes equations. In Section~\ref{section: NumericalDiscretization} we provide both space and time discretizations for the two mentioned FOM, which is essential to get high-fidelity snapshots. In Section \ref{section: NumericalResults} we show the numerical results related to the KPOD-NN technique on the two test cases. Finally, we draw conclusions and possible future developments in Section \ref{section: Conclusions}.

\section{Reduced basis methods for parametrized PDEs using neural networks}
\label{section: KPODNNMethod}

We describe in this section the proposed method, considering for the sake of generality the solution of a nonlinear, time-dependent parametrized PDE. Let us denote by $\Omega \subset \mathbb{R}^d$, $d \geq 1$ a given domain, and by $\mathcal{P} \subset \mathbb{R}^m$ the parameter set, with $m \geq 1$. For each $\boldsymbol{\mu} \in \mathcal{P}$, $\boldsymbol{x} \in \Omega$, $t \in [0, T]$, let us denote by  $\boldsymbol{u}(\boldsymbol{x}, t; \boldsymbol{\mu}): \mathbb{R}^{m+d+1} \xrightarrow{} \mathbb{R}^l$ the scalar ($l=1$) or vector ($l>1$) field, solution of the following problem:
\begin{equation}
\left\{
\begin{array}{ll}
\dfrac{\partial \boldsymbol{u}}{\partial t} + \boldsymbol{A}(t; \boldsymbol{\mu}) \boldsymbol{u} + \boldsymbol{N}(\boldsymbol{u}, t; \boldsymbol{\mu})  = \boldsymbol{f}(t; \boldsymbol{\mu}) & \quad \text{in} \;\; \Omega \times (0, T) \smallskip \\
+ \mathrm{boundary} \ \mathrm{conditions}  & \quad \text{in}  \;\; \partial \Omega \times (0, T) \smallskip \\
+  \mathrm{initial}  \ \mathrm{conditions} & \quad \text{in}  \;\; \Omega \times \{t = 0\}.
 \end{array}
 \right.
\label{eqn: BVP}
\end{equation}
Problem \eqref{eqn: BVP} can be either linear (if  $\boldsymbol{N}(\boldsymbol{u}, t; \boldsymbol{\mu}) = \boldsymbol{0}$) or nonlinear (if $\boldsymbol{N}(\boldsymbol{u}, t; \boldsymbol{\mu}) \neq \boldsymbol{0}$); here $\boldsymbol{f}$ acts as a forcing term.

Our final goal is the efficient numerical approximation of problem \eqref{eqn: BVP} for different values of the parameters vector $\boldsymbol{\mu}$. Our starting point is a full order model (FOM) obtained by introducing the  semi-discrete space approximation of \eqref{eqn: BVP}, of the following form:
\begin{equation}
\left\{
\begin{array}{ll}
\dfrac{\partial \boldsymbol{u}_h}{\partial t} + \boldsymbol{A}_h(t; \boldsymbol{\mu}) \boldsymbol{u}_h + \boldsymbol{N}_h(\boldsymbol{u}_h, t; \boldsymbol{\mu})  = \boldsymbol{f}_h(t; \boldsymbol{\mu}) & \quad \text{for} \;\; t \in (0, T) \smallskip \\
+  \mathrm{initial}  \ \mathrm{conditions} & \quad \text{for}  \;\; t = 0.
 \end{array}
 \right.
\label{eqn: BVPDiscrete}
\end{equation}
This latter problem is generated by introducing a suitable computational grid over the domain $\Omega$, and an (e.g., finite element) approximation $\boldsymbol{u}_h(\boldsymbol{x}, t; \boldsymbol{\mu})$ of $\boldsymbol{u}(\boldsymbol{x}, t; \boldsymbol{\mu})$ depending on a set of $N_h$ degrees of freedom (DOFs); for this reason, with a slight abuse of notation, we denote by $\boldsymbol{u}_h( t; \boldsymbol{\mu}) \in \mathbb{R}^{N_h}$ the DOFs vector corresponding to the solution; $\boldsymbol{A}_h (t; \boldsymbol{\mu}) \in \mathbb{R}^{N_h \times N_h}$ defines the stiffness matrix obtained by the discretization of the (linear) operator $\boldsymbol{A}$, $\boldsymbol{N}_h(\boldsymbol{u}_h, t; \boldsymbol{\mu}) \in \mathbb{R}^{N_h}$ is the vector obtained by discretizing the nonlinear operator ${\bf N}$, and $\boldsymbol{f}_h(t; \boldsymbol{\mu}) \in \mathbb{R}^{N_h}$ corresponds to the discrete source term.
We then consider a partition of the time interval $[0, T]$ into $N_T+1$ equally spaced values $\{t_0 = 0, t_1, \ldots, t_j, ..., t_{N_T} = T\}$, and approximate the time-derivative using suitable finite difference schemes, finally yielding the fully discretized problem, which provides our FOM. For each $\boldsymbol{\mu} \in \mathcal{P}$, the FOM solution at time $t_j$ is then denoted by $\boldsymbol{u}_h^j(\boldsymbol{\mu}) \in \mathbb{R}^{N_h}$, $j \in \{ 0, 1, \ldots, N_T \}$.

Several reduced order models (ROMs) for parametrized PDEs are built by sampling the parameter space and computing a set of {\em snapshots}, that is, FOM solutions for sampled values $\mathcal{P}_h = \{ \boldsymbol{\mu}_1, \boldsymbol{\mu}_2, ..., \boldsymbol{\mu}_i, ..., \boldsymbol{\mu}_m \} \subset \mathcal{P}$ of the parameters. For the case at hand, the  snapshots matrix collecting the computed FOM solutions includes both time and parameter dependencies, and reads:
\begin{equation}
    \mathbb{S} = [ \boldsymbol{u}_h^0(\boldsymbol{\mu}_1) \, | \, \boldsymbol{u}_h^1(\boldsymbol{\mu}_1) \, | \ldots \, | \, \boldsymbol{u}_h^{N_T}(\boldsymbol{\mu}_1) \, | \, \ldots \, | \,\boldsymbol{u}_h^0(\boldsymbol{\mu}_m) \, | \, \boldsymbol{u}_h^1(\boldsymbol{\mu}_m) \, |  \, \ldots \,  | \, \boldsymbol{u}_h^{N_T}(\boldsymbol{\mu}_{m})] \in \mathbb{R}^{N_h \times N_s},
\label{eqn: SnapshotMatrix}
\end{equation}
with $N_s = m \times (N_T + 1)$. For the sake of simplicity, we suppose that $\mathbb{S}$ is built using all the snapshots in time generated by solving the FOM with a fixed time-step $\Delta t$, but eventually only a subset of those solutions can be considered, and, to make the notation lighter, we denote the snapshots matrix by:
\[
\mathbb{S} = [ {\bf s}_1 \, | \, \ldots \, | \, {\bf s}_{N_s}].
\]

To build the proposed ROM strategy, we start from a projection-based framework, exploiting the reduced basis method for parametrized PDEs. We denote by  $\{\boldsymbol{q}_1, \boldsymbol{q}_2, \ldots, \boldsymbol{q}_n \} \subset \mathbb{R}^{N_h}$, with $n \ll N_h$, a set of reduced basis functions spanning a low-dimensional subspace $V_n$, in which the reduced order approximation of the problem solution is sought. Such a subspace
 can be built in two different ways, either linearly, for instance through proper orthogonal decomposition (POD) \cite{Hesthaven1, Hesthaven2}, or nonlinearly (see, e.g., \cite{Fresca1, Fresca2, Lee}). We recall that, in the case of POD, the basis functions can be equivalently obtained either as left singular vectors of the snapshots matrix $\mathbb{S}$, or as eigenvectors (up to a normalization factor) of the  Gram matrix $\mathbb{C} = \mathbb{S}^{\top} \mathbb{S}$.

\begin{figure}[t!]
  \centering
  \includegraphics[keepaspectratio, width=0.7\textwidth]{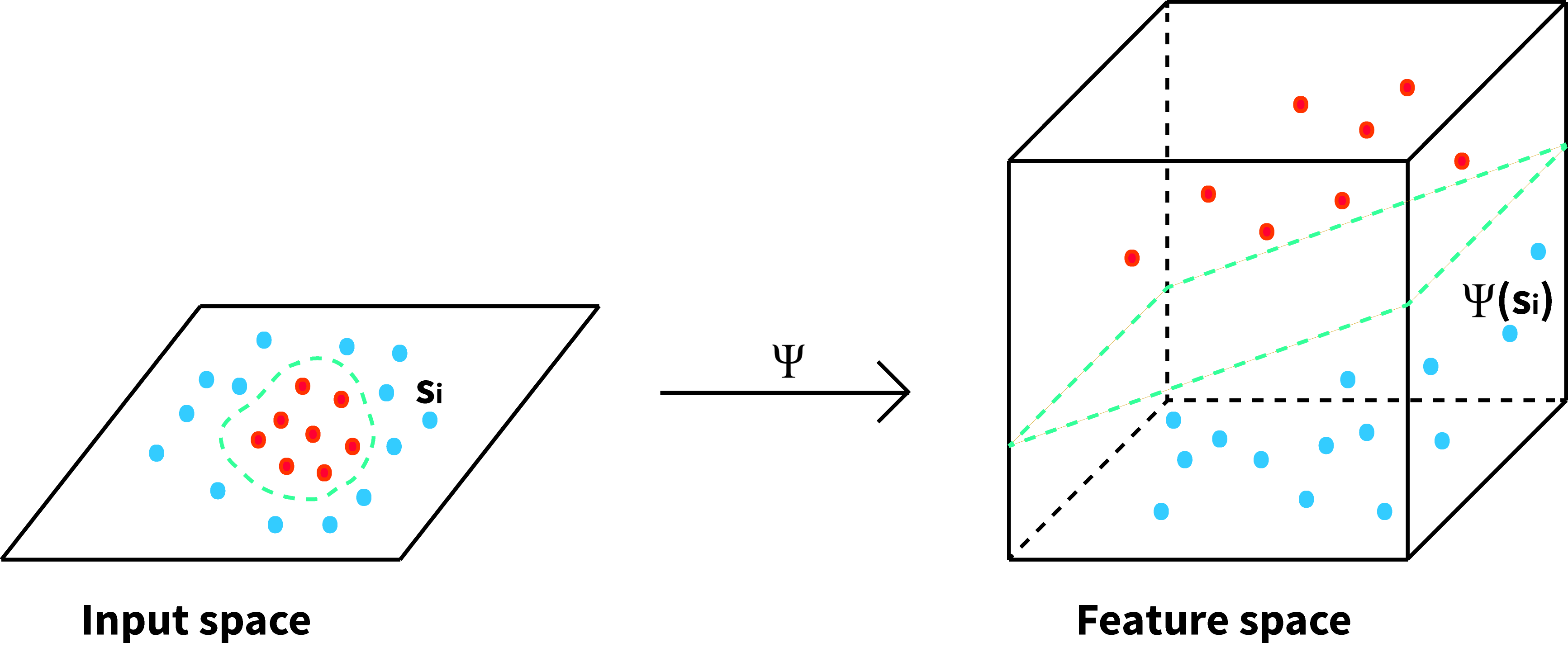}
  \caption{Example of application of the Cover's theorem. Linear inseparable data become separable if projected upon a higher dimensional space through a map $\boldsymbol{\Psi}(\cdot)$.}
  \label{fig: Cover}
\end{figure}

 In this paper, we propose to rely on KPOD \cite{Scholkopf}, which is the nonlinear counterpart of proper orthogonal decomposition (POD), as an efficient and effective alternative to POD.
Compared to POD, KPOD has the additional advantage to be able to deal with otherwise linearly inseparable data. In particular, KPOD projects, through a nonlinear map, the FOM snapshots to a high dimensional space, where we are likely to obtain linear separability (according to Cover's theorem \cite{Cover}), as depicted in Figure \ref{fig: Cover}.

In principle, we should employ an arbitrary transformation $\boldsymbol{\Psi} : \mathbb{R}^{N_h} \rightarrow \mathbb{R}^{N_Z}$ for some very large dimension $N_Z \gg N_h$. Such a transformation is meant to flatten the nonlinear solution manifold $\mathcal{S}_h = \{ \boldsymbol{u}_h( t; \boldsymbol{\mu}), \ t \in [0,T], \ \boldsymbol{\mu} \in \mathcal{P}\} \subset \mathbb{R}^{N_h}$ where the FOM solution is sought. In other words, the manifold obtained as  $\mathcal{Z}_h  = \boldsymbol{\Psi} (\mathcal{S}_h ) = \{  \boldsymbol{U}_h ( t; \boldsymbol{\mu}) = \boldsymbol{\Psi} (\boldsymbol{u}_h( t; \boldsymbol{\mu})) , \ t \in [0,T], \ \boldsymbol{\mu} \in \mathcal{P}\} \subset \mathbb{R}^{N_Z}$ should be more readily  approximated through linear subspaces. To this aim, we would need to apply POD to the $N_Z \times N_s$ matrix containing the transformed snapshots:
\[
\mathbb{S}_{\boldsymbol{\Psi}} = [ \boldsymbol{\Psi} ({\bf s}_1)  \,  | \ldots \, | \,   \boldsymbol{\Psi} ({\bf s}_{N_s} ) ] \in \mathbb{R}^{N_Z \times N_s}
\]
However, we could also apply POD to the Gram matrix $\tilde{\mathbb{C}}_{\boldsymbol{\Psi}} = \mathbb{S}_{\boldsymbol{\Psi}}^{\top} \mathbb{S}_{\boldsymbol{\Psi}}$, which is of size $N_s \times N_s$ as $\mathbb{C} = \mathbb{S}^{\top} \mathbb{S}$. In the context of KPOD, the so-called {\em kernel trick} is applied, i.e. the matrix $\mathbb{S}_{\boldsymbol{\Psi}}$ is directly defined by introducing a bivariate symmetric form $\kappa: \mathbb{R}^{N_h} \times \mathbb{R}^{N_h} \rightarrow \mathbb{R}$, also referred to as {\em kernel function}, rather than the map $\boldsymbol{\Psi}$, to avoid the computation of the components of the matrix  $\tilde{\mathbb{C}}_{\boldsymbol{\Psi}}$ according to the definition, that is:
\[
(\tilde{\mathbb{C}}_{\boldsymbol{\Psi}})_{ij} =  \sum_{k=1}^{N_Z} (\mathbb{S}_{\boldsymbol{\Psi}})_{ki} (\mathbb{S}_{\boldsymbol{\Psi}})_{kj} =
 ( \boldsymbol{\Psi} ({\bf s}_i))^{\top}   ( \boldsymbol{\Psi} ({\bf s}_j)) \qquad i,j = 1,\ldots,N_s.
\]
Thanks to the kernel function, we can instead define:
\[
(\tilde{\mathbb{C}})_{ij} =  \kappa ({\bf s}_i, {\bf s}_j) \qquad i,j = 1,\ldots,N_s,
\]
which takes the name of  kernel (similarity) matrix.
A common choice is to use the squared exponential (or radial basis function) kernel function, which yields the following kernel matrix $\mathbb{K}$:
\begin{equation}
\mathbb{K}({\bf s}_i, {\bf s}_j) = \mathrm{exp} \left( - \gamma ||{\bf s}_i - {\bf s}_j ||_2^2 \right) \in \mathbb{R}^{N_{s} \times N_{s}} \;\;\; \;\;\; i,j = 1, \ldots, N_s.
\label{eqn: KernelMatrix}
\end{equation}
Here $\gamma > 0$ is a hyperparameter of the model, whereas $\boldsymbol{s}_i$ and $\boldsymbol{s}_j$ are two columns of the snapshots matrix $\mathbb{S}$, i.e. two different samples in the time domain and in the parameters space.

One of the advantages of using KPOD is that it does not work in the high-dimensional latent space directly. Indeed, it computes the projections of the FOM solutions onto the principal components, but not the principal components themselves (as done instead in POD, or PCA). To evaluate these projections, first, we compute the eigenvectors and eigenvalues of \eqref{eqn: KernelMatrix}:
\begin{equation}
    \mathbb{K} \boldsymbol{w}_k = \lambda_k \boldsymbol{w}_k \;\;\; \;\;\; k = 1, \ldots, N_s,
    \label{eqn: eig}
\end{equation}
where we denote by $\mathbb{W} = \left[\boldsymbol{w}_1 \, | \,  \ldots \, | \, \boldsymbol{w}_{N_{s}} \right] \in \mathbb{R}^{N_{s} \times N_{s}}$ the  matrix collecting the eigenvectors of $\mathbb{K}$, and by $\Lambda = \text{diag}(\lambda_1, \ldots, \lambda_{N_{s}}) \in \mathrm{R}^{N_{s} \times N_{s}}$ the matrix of the corresponding eigenvalues. Then, by analogy with POD, we compute and collect the projected vectors:

\begin{equation}
  \boldsymbol{p}_k = \frac{1}{\sigma_k} \mathbb{S}  \boldsymbol{w}_k \qquad k=1,\ldots,N_s
\end{equation}
in the matrix $\mathbb{P} = \left[ \boldsymbol{p}_1 \,  | \, \ldots \, | \, \boldsymbol{p}_{N_s} \right] \in \mathbb{R}^{N_h \times N_s}$, with $\sigma_k = \sqrt{\lambda_k}$.
In POD, $\mathbb{P}$ would be an orthogonal matrix.
Conversely, in KPOD, we finally have to perform a reduced QR factorization of $\mathbb{P}$ to obtain an orthonormal basis \cite{Quarteroni1}:
\begin{equation}
    \left[\mathbb{Q}, \tilde{\mathbb{R}}\right] = \text{QR}(\mathbb{P}),
\label{eqn: reducedQRfactorization}
\end{equation}
where $\mathbb{Q} \in \mathbb{R}^{N_h \times N_{hs}}$ is an orthogonal matrix, whereas $\tilde{\mathbb{R}} \in \mathbb{R}^{N_{hs} \times N_s}$ is an upper triangular matrix, being $N_{hs} = \text{min}(N_h, N_s)$.

We select $n$ columns of $\mathbb{Q}$ to get the (orthonormal) reduced basis $\mathbb{V} = [\boldsymbol{q}_1 | \boldsymbol{q}_2 | ... | \boldsymbol{q}_n] \in \mathbb{R}^{N_h \times n}$.
The dimension $n$ where the truncation occurs is determined as follows:
\begin{equation}
\dfrac{\sum_{k=n+1}^{N_s} \sigma_k^2}{\sum_{k=1}^{N_s} \sigma_k^2} \leq \hat{\varepsilon},
\label{eqn: StoppingCriterion}
\end{equation}
being $\hat{\varepsilon} > 0$ a given tolerance and $\sigma_k^2$ the square of the $k^{th}$ singular value. We will always use $\hat{\varepsilon} = 10^{-12}$ in the numerical results.

The process is summarized in Algorithm 1. So far, we have just described the forward mapping, i.e. the nonlinear dimensionality reduction process. Due to the nonlinear nature of this technique, it is not trivial to map an element from the reduced space (with dimension $n$) towards its {\em pre-image} in the input space of the high-fidelity solutions (with dimension $N_h$). We refer to this mathematical problem as KPOD backward mapping. %

\begin{figure}[t!]
  \centering
  \includegraphics[keepaspectratio, width=0.75\textwidth]{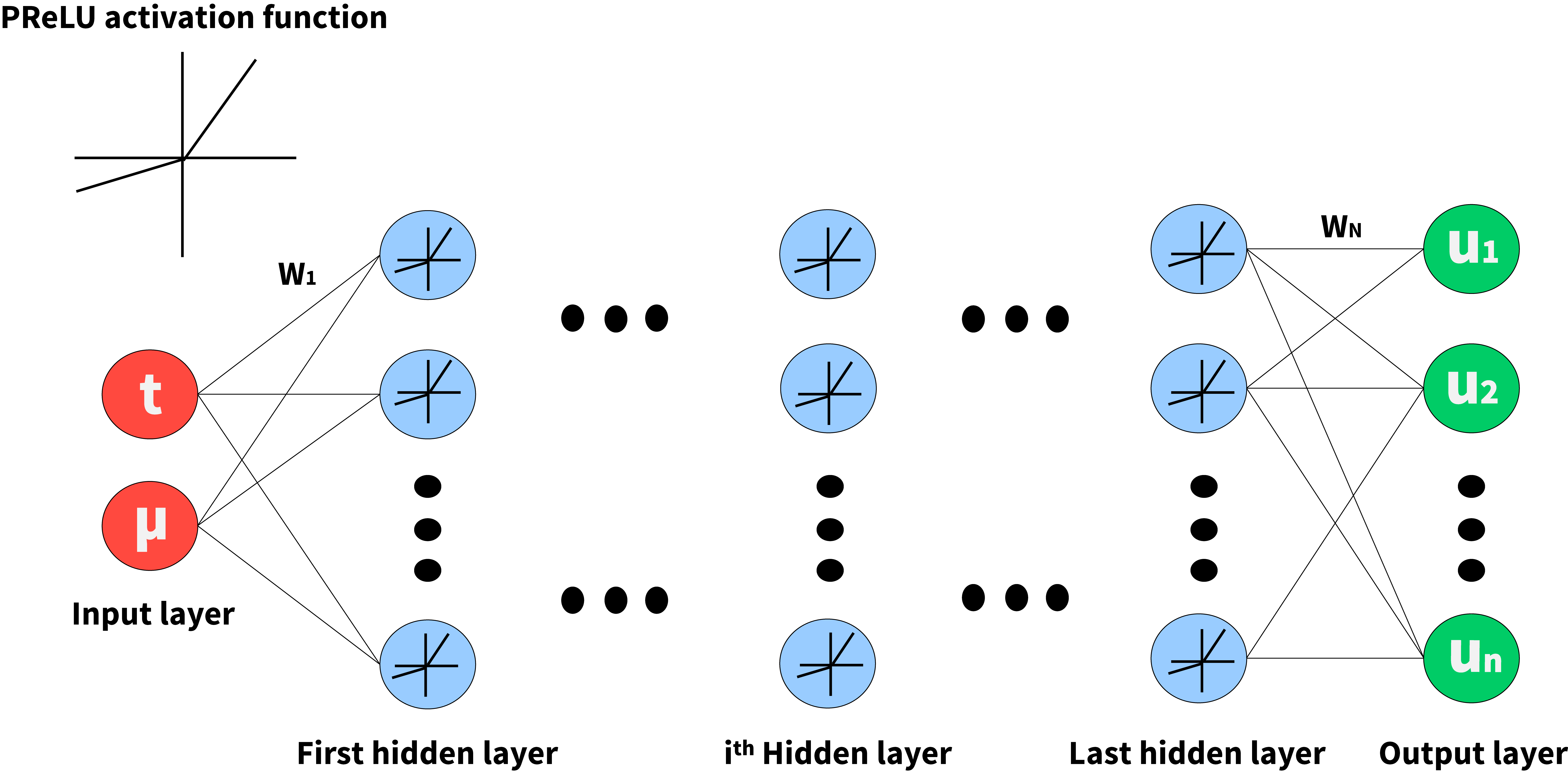}
  \caption{Example of fully connected feedforward neural network with PReLU activation functions, several hidden layers and an input layer made by time and one single parameter $\mu$ ($m=1$). $\boldsymbol{W}_i$ indicates weights of the connections between the $i^{th}$ layer and the following one.}
  \label{fig: NNArchitecture}
  \vspace{3mm}
\end{figure}

Once $\mathbb{V}$ is built, we obtain the reduced coefficients by means of the linear projection $\mathbb{V}^T \boldsymbol{u}_h(t; \boldsymbol{\mu})$ onto $V_n$.
We introduce now a function $\boldsymbol{\Phi}$ which represents the connection between the space of parameters, wherein we assimilate the time independent variable to a parameter, and the projected coefficients, i.e:
\begin{equation}
    \boldsymbol{\Phi}: \mathcal{T} \times \mathcal{P} \subset \mathbb{R}^{m+1} \xrightarrow{} \mathbb{R}^n \;\; \mathrm{s.t.} \;\; \boldsymbol{\Phi}: (t; \boldsymbol{\mu}) \xrightarrow{} \boldsymbol{\mathbb{V}}^T \boldsymbol{u}_h(t; \boldsymbol{\mu}).
    \label{eqn: ParametricMap}
\end{equation}
Similarly to what is done in \cite{Hesthaven2} for the POD-NN method, we propose to rely on NN to determine the reduced order approximation for any new parameter instance.
This NN is able to perform a nonlinear regression for the KPOD-NN method, to learn an approximation $\boldsymbol{\Phi}_{NN}$ of the map $\boldsymbol{\Phi}$. The evaluation of $\boldsymbol{\Phi}_{NN}$ for a certain $(t, \boldsymbol{\mu})$ provides a reduced solution:
\begin{equation}
  \boldsymbol{u}_n^{NN}(t; \boldsymbol{\mu}) = \boldsymbol{\Phi}_{NN}(t; \boldsymbol{\mu}),
\end{equation}
and consequently we come back to the high-fidelity dimension through:
\begin{equation}
  \boldsymbol{u}_h^{NN}(t; \boldsymbol{\mu}) = \mathbb{V} \boldsymbol{\Phi}_{NN}(t; \boldsymbol{\mu}) = \mathbb{V} \boldsymbol{u}_n^{NN}(t; \boldsymbol{\mu}).
\end{equation}
To the best of our knowledge, a novelty is that we consider a NN where the number of hidden layers scales as $\ceil*{\text{log}(n)}$ and the size of each hidden layer depends on the dimension $n$ of the reduced basis. We also have a first (input) layer with $N_{I} = m+1$ units (for the $(t; \boldsymbol{\mu})$ values) and a final (output) layer of size $N_{O} = n$ (for each single component of the reduced solution $\mathbb{V}^T \boldsymbol{u}_h(t; \boldsymbol{\mu})$). Our NN is fully connected, i.e. each neuron of a specific layer is connected to all the neurons of the next layer only \cite{Hesthaven1}. The strength of these interactions is defined by weights $\boldsymbol{W}_i$ and biases \cite{Goodfellow} (as shown in Figure \ref{fig: NNArchitecture}), whose values change during the training process. In this way the NN provides a suitable lower dimensional representation of the specific FOM that we want to approximate. We use a Parametric Rectified Linear Unit (PReLU) as activation function of the hidden layers. This function is known to be more flexible than ReLU and LeakyReLU while leading better NN approximations \cite{Pedamonti}. PreLU reads:
\begin{equation}
\mathrm{PReLU}(x) = \begin{cases}
                      \alpha x & $in$ \; x < 0, \\
                      x & $in$ \; x \geq 0,
                    \end{cases}
\end{equation}
where $\alpha > 0$ is automatically optimized during the training phase.

The architecture is synthesized in Table \ref{table: NNArchitecture}. A schematic view of the NN is sketched in Figure \ref{fig: NNArchitecture}.

\begin{algorithm}[H]
\SetAlgoLined
 \textbf{function} $\mathbb{V}$=KPOD($\mathbb{S}$, $\hat{\varepsilon}$, $\gamma$)
 \begin{algorithmic}
   \Statex Assemble $\mathbb{K}({\bf s}_i, {\bf s}_j) = \mathrm{exp} \left( - \gamma ||{\bf s}_i - {\bf s}_j ||_2^2 \right) \;\;\; \text{for} \;\;\; i,j = 1, ..., N_s$.
   \Statex Perform an eigenvalue-eigenvector decomposition: $\mathbb{K} \boldsymbol{w}_k = \lambda_k \boldsymbol{w}_k \;\;\; \text{for} \;\;\; k = 1, ..., N_s$.
   \Statex Compute $\mathbb{P} = \left[ \boldsymbol{p}_1 | \boldsymbol{p}_2 | ... | \boldsymbol{p}_{N_s} \right] \in \mathbb{R}^{N_h \times N_s}$, where $\boldsymbol{p}_k = \frac{1}{\sqrt{\lambda_k}} \mathbb{S}  \boldsymbol{w}_k \;\;\; \text{for} \;\;\;  k=1,\ldots,N_s$.
   \Statex Reduced QR factorization of matrix $\mathbb{P}$: $\left[\mathbb{Q}, \tilde{\mathbb{R}}\right] = \text{QR}(\mathbb{P})$.
   \Statex Select $n$ $\boldsymbol{q}_k$ vectors according to \eqref{eqn: StoppingCriterion}.
   \Statex Assemble $\mathbb{V} = [\boldsymbol{q}_1 | \boldsymbol{q}_2 | ... | \boldsymbol{q}_n]$.
 \end{algorithmic}
 \textbf{end function}
 \label{algorithm: KPOD}
 \caption{KPOD as nonlinear dimensionality reduction technique for ROM.}
\end{algorithm}

\begin{table}[h!]
\centering
\begin{tabular}{|c c|}
    \hline
        Layer & Number of neurons \\ [0.5ex]
    \hline
        Dense + PReLU activation (input) & $m+1$ \\
        Dense + PReLU activation (hidden) 1 & $n$ \\
        Dense + PReLU activation (hidden) 2 & $n$ \\
        ... & ... \\
        Dense + PReLU activation (hidden) $\ceil*{\text{log}(n)}$ & $n$ \\
        Dense (output) & $n$ \\
    \hline
\end{tabular}
\vspace{10pt}
\caption{NN architecture. $n$ is the dimension of the reduced basis, whereas $m$ is the length of the parameter vector. Glorot uniform initializer \cite{Glorot2010} is used to set the initial values of the weights of all layers.}
\label{table: NNArchitecture}
\end{table}

For the sake of NN construction, in supervised learning, let $\mathcal{D}_{b} = \{\boldsymbol{x}_i, \boldsymbol{y}_i\}_{i = 1}^{N_{b}}$ be a building set with input-output pairs, which is split in different subsets: a training (equivalently learning) dataset $\mathcal{D}_{tr} = \{\boldsymbol{x}_i, \boldsymbol{y}_i\}_{i = 1}^{N_{tr}}$ and a validation one $\mathcal{D}_{v} = \{\boldsymbol{x}_i, \boldsymbol{y}_i\}_{i = 1}^{N_{v}}$, with $N_{tr} + N_{v} = N_{b}$. Inputs $\boldsymbol{x}_i$ are sampled from the space of parameters, whereas outputs $\boldsymbol{y}_i$ represent the corresponding reduced coefficients. We apply the so called $K$-fold cross-validation \cite{Hastie}, where the building set is divided into $K$ equal parts (here we use $K=5$). In this respect, the network is trained $K$ times, where $K-1$ subdivisions of the building set act as learning set, whereas the last remaining subdivision is used to validate the model for a specific choice of the hyperparameters. %
Data located inside the building set are shuffled at the beginning of the $K$-fold procedure so that there is a high probability for the network to learn the reduced numerical solution along the entire time span $[0, T]$.

Thanks to $K$-fold cross-validation, we are able to perform both the training phase and the model selection.
In particular, we investigate the effects of changing the number of layers and the growth/decay of the number of neurons of the layers in the network (e.g. linear, parabolic and hyperbolic). The selection of the best model among the different analyzed architectures is made through the comparison of their generalization error $E_{gen}$:
\begin{equation}
E_{gen} = \dfrac{1}{N_v} \sum_{k=1}^{K} \dfrac{\sqrt{\sum_{i=1}^{N_v} (\boldsymbol{\hat{y}}_i^k - \boldsymbol{y}_i)^2}}{\sqrt{\sum_{i=1}^{N_v} \boldsymbol{y}_i^2}},
\end{equation}
where $\boldsymbol{\hat{y}}_i^k$ is the prediction of the neural network on fold $k$ associated with the output observation $\boldsymbol{y}_i$. We look for the network with the smallest generalization error. We trained different neural networks using input-output pairs coming from the benchmark problems that will be introduced in Section \ref{section: MathematicalModeling}. As the generalization error was not significantly affected by the growth/decay strategy used to determine the number of neurons in a single layer, we selected a constant number of neurons per layer. The selected architecture is reported in Table \ref{table: NNArchitecture}, which provided  - among the tested architectures - the lowest generalization error. In order to adapt the network complexity to the dimension $n$, we scaled the number of layers as $\ceil*{\text{log}(n)}$.

Once the architecture is fixed and a NN is built and trained, we test its performances, feeding the network with a new dataset (called test set) $\mathcal{D}_{test} = \{\boldsymbol{x}_i, \boldsymbol{y}_i\}_{i = 1}^{N_{test}}$, which is made by unseen observations. In this way, we assess the learning capability of the network. All the abovementioned operations are implemented using the Keras Python Deep Learning library \cite{Chollet}. Our loss function is the discrete relative $L^2$ norm:
\begin{equation}
\mathcal{C}_j = \dfrac{\sqrt{\sum_{i=1}^{N_j} (\boldsymbol{\hat{y}}_i - \boldsymbol{y}_i)^2}}{\sqrt{\sum_{i=1}^{N_j} \boldsymbol{y}_i^2}} \;\;\;\; \mathrm{with} \;\;\;\; j \in \{tr, v\},
\label{eqn: LossFunction}
\end{equation}
The network is trained over different steps, which are called epochs, considering each time a small amount of observations of the training set, which are referred to as batches; here we consider $N_{batch} = 10$.
A small value of batch size is known to avoid sharp local minima in the optimization process, improving the accuracy of the network \cite{Keskar}. We use the Adam stochastic optimizer \cite{Kingma} with AMSGrad variant, $\beta_1 = 0.9$ and $\beta_2 = 0.999$, to find the optimal values of weights and bias. The initial learning rate $lr$ is set to $0.01$ for the wave equation and to $0.1$ for the lid-driven cavity problem. This lets the stochastic gradient descent method better explore the landscape of local minima during the very first epochs of training.
We consider a regularization of the weights so that the network does not overfit and is able to explore more possibilities in the training process. The final formulation of the loss function is reported here:
\begin{equation}
\mathcal{C}_j = \dfrac{\sqrt{\sum_{i=1}^{N_j} (\boldsymbol{\hat{y}}_i - \boldsymbol{y}_i)^2}}{\sqrt{\sum_{i=1}^{N_j} \boldsymbol{y}_i^2}} + \Theta \sum_{l=1}^{N_{l}} \sum_{o} \sum_{k} w_{l, o, k}^2 \;\;\;\; \mathrm{with} \;\;\;\; j \in \{tr, v\},
\end{equation}
where $\Theta = 0.01$ acts as regularization parameter, $w_{l, o, k}$ is the $o^{th}$ weight of the $l^{th}$ layer connected to the $k^{th}$ neuron of the following layer of the NN.
The entire implementation of the offline and online stages of the KPOD-NN method is described in Algorithm 2. This method is non-intrusive, because it does not need to access system operators to perform the projection and the online evaluation of the feedforward neural network is independent of the high-fidelity numerical scheme. For the sake of simplicity, we present an approach in which the size of the reduced space plays a role only in the definition of the NN architecture. The underlying idea is that the higher the dimension of the reduced space the higher should be the complexity of the NN. Still, a dependence from the dimension $n$ can be eventually introduced in other parameters of the NN, such as the regularization term $\Theta$, the configuration of the optimizer, the total number of epochs and the batch size.

\begin{algorithm}[H]
\SetAlgoLined
 \textbf{function} $[\mathbb{V}, \boldsymbol{\Phi}_{NN}]$=KPOD-NN\_OFFLINE($\mathcal{P}$, $\Omega$, $T$, $m$, $N_T$, $\hat{\varepsilon}$, $\gamma$)
 \begin{algorithmic}
   \Statex Generate $[t_0 = 0, t_1, ..., t_j, ..., t_{N_T} = T]$ and $\mathcal{P}_h = \{ \boldsymbol{\mu}_1, \boldsymbol{\mu}_2, ..., \boldsymbol{\mu}_i, ..., \boldsymbol{\mu}_m \} \subset \mathcal{P}$.
   \Statex Compute high-fidelity numerical solutions $\boldsymbol{u}_h^j(\boldsymbol{\mu}_i)$, $j = 0, ..., N_T$, $i = 1, ..., m$.
   \Statex Set $\mathbb{S} = [ \boldsymbol{u}_h^0(\boldsymbol{\mu}_1), \boldsymbol{u}_h^1(\boldsymbol{\mu}_1), ..., \boldsymbol{u}_h^{N_T}(\boldsymbol{\mu}_1), \boldsymbol{u}_h^0(\boldsymbol{\mu}_2), ..., \boldsymbol{u}_h^{N_T}(\boldsymbol{\mu}_{m})]$.
   \Statex $\mathbb{V}$=KPOD($\mathbb{S}$, $\hat{\varepsilon}$, $\gamma$).
   \Statex Build and train the NN for the calculation of $\boldsymbol{\Phi}_{NN}$.
 \end{algorithmic}
 \textbf{end function}

 \textbf{function} $\mathbb{V}$=KPOD-NN\_ONLINE($(t_{new}; \boldsymbol{\mu}_{new})$, $\mathbb{V}$, $\boldsymbol{\Phi}_{NN}$)
 \begin{algorithmic}
   \Statex Evaluate the output $\boldsymbol{u}_n^{NN}(t_{new}; \boldsymbol{\mu}_{new}) = \boldsymbol{\Phi}_{NN}(t_{new}; \boldsymbol{\mu}_{new})$ of the trained NN for the input $(t_{new}; \boldsymbol{\mu}_{new})$.
   \vspace{2pt}
   \Statex Map the solution to the high-fidelity space as $\boldsymbol{u}_h^{NN}(t_{new}; \boldsymbol{\mu}_{new}) = \mathbb{V} \boldsymbol{\Phi}_{NN}(t_{new}; \boldsymbol{\mu}_{new})$.
 \end{algorithmic}
 \textbf{end function}
 \caption{KPOD-NN ROM method for unsteady PDEs.}
\end{algorithm}

\section{Benchmark problems}
\label{section: MathematicalModeling}
We apply the method described in Section \ref{section: KPODNNMethod} to the following two differential problems, namely the wave equation and the Navier-Stokes equations. The first problem is set in 1D, whereas the second one will be solved in 2D.
\subsection{Wave equation}
The wave equation is a second-order hyperbolic linear PDE and reads:
\begin{equation}
\begin{cases}
\dfrac{\partial^2 u}{\partial t^2} = c^2 \dfrac{\partial^2 u}{\partial x^2} & $in$ \; \Omega \times (0, T), \\
u(x, 0) = A_0 \mathrm{exp}^{- \frac{(x - x_0)^2}{2 \sigma^2} } & $in$ \; \Omega \times \{0\}, \\
u(0, t) = 0 & t \in (0, T), \\
u(L, t) = 0 & t \in (0, T),
\label{eqn: Wave}
\end{cases}
\end{equation}
with $\Omega = (0, L)$. $A_0 \in [0.5, 1]$, $x_0 \in [ \frac{L}{3}, \frac{2L}{3}]$, $\sigma \in [0.5, 1]$ and $c \in \mathbb{R}$, with $L=4 \pi$ and $T = 52$. We assign homogeneous Dirichlet boundary conditions at $x=0$ and $x=L$. The initial condition, given by a Gaussian pulse, triggers the propagation of two waves in opposite directions. We consider intensity $A_0$, position $x_0$ and variability $\sigma$ of the initial pulse as input parameters of the neural network, i.e. $u = u(x, t; A_0, x_0, \sigma)$.

\subsection{Navier-Stokes equations}

The Navier-Stokes equations model the flow of a viscous fluid, either incompressible or compressible, in a certain domain $\Omega \subset \mathbb{R}^2$ and in the time interval $(0, T)$. The unknowns in primitive variables are velocity and pressure, i.e $(\boldsymbol{u}, p)$. The corresponding dimensionless strong form reads \cite{Tagliabue}:
\begin{equation}
\begin{cases}
\dfrac{\partial \boldsymbol{u}}{\partial t} + (\boldsymbol{u} \cdot \nabla) \boldsymbol{u} + \nabla p -2 \nabla \cdot \left( \dfrac{1}{\mathbb{R}e} D(\boldsymbol{u}) \right) = \boldsymbol{f} & $in$ \; \Omega \times (0, T), \\
\nabla \cdot \boldsymbol{u} = 0 & $in$ \; \Omega \times (0, T), \\
\boldsymbol{u}(0) = \boldsymbol{u}_0 & $in$ \; \Omega \times \{0\}, \\
\boldsymbol{\sigma}_f \boldsymbol{n} = \boldsymbol{h} & $on$ \; \Gamma_N \times (0, T), \\
\boldsymbol{u} = \boldsymbol{g} & $on$ \; \Gamma_D \times (0, T),
\label{eqn: NavierStokesPrimitive}
\end{cases}
\end{equation}
where $\mathbb{R}e$ is the Reynolds number, $D(\boldsymbol{u}) = \tfrac{1}{2} (\nabla \boldsymbol{u} + \nabla \boldsymbol{u}^T)$ is the strain tensor, and $\boldsymbol{\sigma}_f = -p I + 2 \mathbb{R}e^{-1} D(\boldsymbol{u})$ is the Cauchy stress tensor, being $I$ the identity tensor. $\boldsymbol{f}: \Omega \times (0, T) \xrightarrow{} \mathbb{R}^2$ indicates the body forces, $\boldsymbol{u}_0: \Omega \xrightarrow{} \mathbb{R}^2$. $\boldsymbol{h}$ and $\boldsymbol{g}$ denote respectively the vector fields for the Neumann boundary condition on $\Gamma_N$ and the Dirichlet boundary condition on $\Gamma_D$, where $\Gamma_N, \Gamma_D \subseteq \partial \Omega$ with $\overline{\Gamma_D \cup \Gamma_N} = \partial \Omega$ and $\mathring{\Gamma}_D \cap \mathring{\Gamma}_N = \emptyset$. $\boldsymbol{n}$ indicates the outward directed unit vector normal to $\Gamma_N$.

The weak formulation of \eqref{eqn: NavierStokesPrimitive} reads \cite{Tagliabue}: given $V_g = \{ \boldsymbol{v} \in [H^1(\Omega)]^2 : \boldsymbol{v}|_{\Gamma_D} = \boldsymbol{g} \}$, $V = \{ \boldsymbol{v} \in [H^1(\Omega)]^2 : \boldsymbol{v}|_{\Gamma_D} = \boldsymbol{0} \}$, $Q = L^2(\Omega)$, $W_g = \{ \boldsymbol{v} \in V_g : \nabla \cdot \boldsymbol{v} = 0 \} \subset V_g$ and $W = \{ \boldsymbol{v} \in V : \nabla \cdot \boldsymbol{v} = 0 \} \subset V$ spaces of divergence-free functions of $V_g$ and $V$ respectively, find $\boldsymbol{u} = \boldsymbol{u}(t) \in W_g \;\; \forall t \in (0, T)$ such that:
\begin{equation}
\begin{cases}
m\left(\dfrac{\partial \boldsymbol{u}}{\partial t}, \boldsymbol{v}\right) + a(\boldsymbol{u}, \boldsymbol{v}) + b(\boldsymbol{v}, p) + c(\boldsymbol{u}, \boldsymbol{u}, \boldsymbol{v}) = F(\boldsymbol{v}) + H(\boldsymbol{v}) & \forall \boldsymbol{v} \in W, \\
b(\boldsymbol{u}, q) = 0 & \forall q \in Q, \\
\boldsymbol{u}(0) = \boldsymbol{u}_0 & $in$ \; \Omega,
\label{eqn: NavierStokesPrimitiveWeak}
\end{cases}
\end{equation}
where $m\left( \dfrac{\partial \boldsymbol{u}}{\partial t}, \boldsymbol{v} \right) = \displaystyle \int_{\Omega} \dfrac{\partial \boldsymbol{u}}{\partial t} \cdot \boldsymbol{v} \, d\Omega$, $a(\boldsymbol{u}, \boldsymbol{v}) = \dfrac{2}{\mathbb{R}e} \displaystyle \int_{\Omega} D(\boldsymbol{u}):D(\boldsymbol{v}) \, d\Omega$, $b(\boldsymbol{v}, p) = \displaystyle - \int_{\Omega} p \nabla \cdot \boldsymbol{v} \, d\Omega$, $c(\boldsymbol{w}, \boldsymbol{u}, \boldsymbol{v}) = \displaystyle \int_{\Omega} ((\boldsymbol{w} \cdot \nabla) \boldsymbol{u}) \cdot \boldsymbol{v} \, d \Omega$, while $F(\boldsymbol{v}) = \displaystyle \int_{\Omega} \boldsymbol{f} \cdot \boldsymbol{v} \, d\Omega$ and $H(\boldsymbol{v}) = \displaystyle \int_{\Gamma_N} \boldsymbol{h} \cdot \boldsymbol{v} \, d\Gamma$.

The Navier-Stokes equations can be rewritten in the streamfunction formulation \cite{Tagliabue}. We introduce the quotient space $X = \{ [\phi] : \phi \in H^2(\Omega) \}$ of scalar functions in $H^2(\Omega)$ that differ up to a constant and, under the hypothesis of simply connected domain, we consider the unique potential $\phi \in X$ such that $\boldsymbol{v} = \boldsymbol{curl} \phi$, where $\boldsymbol{curl}(\cdot)$ operator is defined as follows:
\begin{equation}
\boldsymbol{curl}: X \xrightarrow{} [H^1(\Omega)]^2 \;\;\;\; \boldsymbol{curl} \phi = \left( \dfrac{\partial \phi}{\partial y}, - \dfrac{\partial \phi}{\partial x} \right)
\label{eqn: curl}
\end{equation}
The streamfunction weak formulation reads: given $\boldsymbol{\Phi}_g = \{\psi \in X: \boldsymbol{curl} \psi |_{\Gamma_D} = \boldsymbol{g} \}$ and $\boldsymbol{\Phi} = \{\psi \in X: \boldsymbol{curl} \psi |_{\Gamma_D} = \boldsymbol{0} \}$, find $\phi = \phi(t) \in \boldsymbol{\Phi}_g \;\; \forall t \in (0, T)$ such that:
\begin{equation}
\begin{cases}
\overline{m}\left(\dfrac{\partial \phi}{\partial t}, \psi \right) + \overline{a}(\phi, \psi) + \overline{c}(\phi, \phi, \psi) = \overline{F}(\psi) + \overline{H}(\psi) & \forall \psi \in \boldsymbol{\Phi}, \\
\phi(0) = \phi_0 & $in$ \; \Omega,
\label{eqn: NavierStokesStreamWeak}
\end{cases}
\end{equation}
where $\overline{m}\left( \dfrac{\partial \phi}{\partial t}, \psi \right) = \displaystyle \int_{\Omega} \boldsymbol{curl} \dfrac{\partial \phi}{\partial t} \cdot \boldsymbol{curl} \psi \, d\Omega$, $\overline{a}(\phi, \psi) = \dfrac{2}{\mathbb{R}e} \displaystyle \int_{\Omega} D(\boldsymbol{curl} \phi):D(\boldsymbol{curl} \psi) \, d\Omega$, $c(\phi, \phi, \psi) = \displaystyle \int_{\Omega} ((\boldsymbol{curl} \phi \cdot \nabla) \boldsymbol{curl} \phi) \cdot \boldsymbol{curl} \psi \, d \Omega$, while $\overline{F}(\psi) = \displaystyle \int_{\Omega} \boldsymbol{f} \cdot \boldsymbol{curl} \psi \, d\Omega$ and $\overline{H}(\psi) = \displaystyle \int_{\Gamma_N} \boldsymbol{h} \cdot \boldsymbol{curl} \psi \, d\Gamma$. The initial condition is chosen in such a way that $\boldsymbol{curl} \phi_0 = \boldsymbol{u}_0$. By construction the velocity field $\boldsymbol{u} = \boldsymbol{curl} \phi$ is divergence free.

\begin{figure}[t!]
	\centering
	\includegraphics[width=0.45\textwidth]{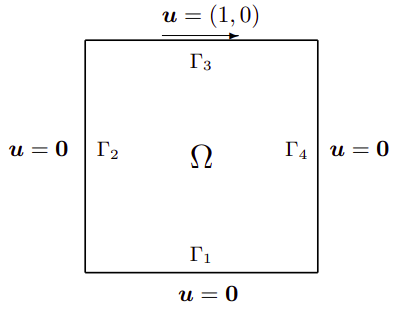}
	\caption{Lid-driven cavity problem: data and settings. $\Omega = (0, 1)^2$, $\Gamma_1 \cup \Gamma_2 \cup \Gamma_3 \cup \Gamma_4 = \partial \Omega$.}
	\label{fig: liddrivencavity}
\end{figure}

We solve the lid-driven cavity benchmark problem on a square domain $\Omega = (0, 1)^2$ as done in \cite{Tagliabue}. In Figure \ref{fig: liddrivencavity}, we depict the boundary conditions assigned to $\partial \Omega$. We consider one single parameter, i.e. the Reynolds number $\mathbb{R}e \in [100, 5000]$, using the streamfunction formulation of the Navier-Stokes equations in 2D. We look for the steady state of the Navier-Stokes equations, i.e. we always take the snapshot at the final time $t=T=30$ only. Indeed, we apply the KPOD-NN reduction on $\phi = \phi(\boldsymbol{x}, T; \mathbb{R}e)$.

\section{Numerical discretization}
\label{section: NumericalDiscretization}

In this section we discretize in space and time problems \eqref{eqn: Wave} and \eqref{eqn: NavierStokesStreamWeak}. %

\subsection{Wave equation}
We use the finite difference method \cite{Quarteroni1} to discretize the wave equation. We consider a partition $[x_0 = 0, x_1, ..., x_N = L]$ of $N + 1$ equally spaced points distributed in the computational domain $\Omega = (0, L)$, and a set $[t_0 = 0, t_1, ..., t_{N_T} = T]$ of $N_T + 1$ equally spaced times in $[0, T]$. We employ second-order centered finite differences to approximate the derivatives of the PDE, leading to the following formulation:
\begin{equation}
\begin{cases}
\dfrac{u_i^{j+1} - 2 u_i^j + u_i^{j-1}}{\Delta t^2} = c^2 \dfrac{u_{i+1}^j - 2 u_i^j + u_{i-1}^j}{\Delta x^2} & \forall i = 1, ..., N-1 \quad \forall j = 0, ..., N_T-1, \\
\displaystyle u_i^0 = A_0 \mathrm{exp}^{- \frac{(x_i - x_0)^2}{2 \sigma^2}} & \forall i = 1, ..., N-1, \\
u_0^j = 0 & \forall j = 0, ..., N_T, \\
u_N^j = 0 & \forall j = 0, ..., N_T,
\label{eqn: WaveFD}
\end{cases}
\end{equation}
where $\Delta t = T/N_T$ and $\Delta x = L/N$. We rearrange the first equation in \eqref{eqn: WaveFD} to obtain $u_i^{j+1}$:
\begin{equation}
u_i^{j+1} = C (u_{i+1}^j + u_{i-1}^j) + 2 (1 - C) u_i^j - u_i^{j-1},
\label{eqn: WaveFirstEquation}
\end{equation}
with $C = c^2 \dfrac{\Delta t^2}{\Delta x^2}$. We remind that the value of $c \dfrac{\Delta t}{\Delta x}$ must be less or equal than 1 to satisfy the Courant-Friedrichs-Lewy condition, which ensures the stability of the numerical scheme \cite{Quarteroni1}.

\subsection{Navier-Stokes equations}

Following the approach proposed in \cite{Tagliabue}, we provide the space discretization of \eqref{eqn: NavierStokesStreamWeak} by means of the NURBS-based Isogeometric Analysis (IGA) Galerkin method. We employ the generalized-$\alpha$ scheme for its time approximation \cite{Chung, Jensen}.

Given an exact representation of the square domain $\Omega$ through B-splines functions, we introduce the finite dimensional space of B-splines in the physical domain \cite{Tagliabue}, say $\mathcal{V}_h$, and we define $\boldsymbol{\Phi}_h = \boldsymbol{\Phi} \cap \mathcal{V}_h$ and $\boldsymbol{\Phi}_{g, h} = \boldsymbol{\Phi}_g \cap \mathcal{V}_h$, being $N_s = dim(\boldsymbol{\Phi})$ the dimension and $\{ \mathcal{R}_i \}_{i=1}^{N_s}$ the set of basis functions. The IGA approximation of \eqref{eqn: NavierStokesStreamWeak} reads: find $\phi_h(t) \in \boldsymbol{\Phi}_{g, h} \;\; \forall t \in (0, T)$ such that:
\begin{equation}
\begin{cases}
\overline{m}\left(\dfrac{\partial \phi_h}{\partial t}, \psi_h \right) + \overline{a}(\phi_h, \psi_h) + \overline{c}(\phi_h, \phi_h, \psi_h) = \overline{F}(\psi_h) + \overline{H}(\psi_h) & \forall \psi_h \in \boldsymbol{\Phi}_h, \\
\phi_h(0) = \phi_{h, 0} & $in$ \; \Omega.
\label{eqn: NavierStokesStreamSpaceDiscretization1}
\end{cases}
\end{equation}

We can rewrite \eqref{eqn: NavierStokesStreamSpaceDiscretization1} in residual form: find $\phi_h(t) \in \boldsymbol{\Phi}_{g, h} \;\; \forall t \in (0, T)$ such that:
\begin{equation}
\begin{cases}
R_h \left( \psi_h, \dfrac{\partial \phi_h}{\partial t}, \phi_h \right) = 0 & \forall \psi_h \in \boldsymbol{\Phi}_h, \\
\phi_h(0) = \phi_{h, 0} & $in$ \; \Omega,
\label{eqn: NavierStokesStreamSpaceDiscretization2}
\end{cases}
\end{equation}
with $R_h \left( \psi_h, \dfrac{\partial \phi_h}{\partial t}, \phi_h \right) = \overline{m}\left(\dfrac{\partial \phi_h}{\partial t}, \psi_h \right) + \overline{a}(\phi_h, \psi_h) + \overline{c}(\phi_h, \phi_h, \psi_h) - \overline{F}(\psi_h) - \overline{H}(\psi_h)$. We define the vector of discrete residuals whose components are the residuals $R_h \left( \cdot, \dfrac{\partial \phi_h}{\partial t}, \phi_h \right)$ evaluated in the NURBS basis functions $\mathcal{R}_i \;\; i=1, ..., N_s$ for the function space $\boldsymbol{\Phi}_h$, i.e. $\boldsymbol{R}\left( \dfrac{\partial \phi_h}{\partial t}, \phi_h \right) = \bigg\{R_h \left( \mathcal{R}_i, \dfrac{\partial \phi_h}{\partial t}, \phi_h \right) \bigg\}_{i=1}^{N_s}$, $\forall t \in (0, T)$. Moreover, we introduce $\boldsymbol{\phi}=\boldsymbol{\phi}(t)=\{ \phi_i \}_{i=1}^{N_s}$ and $\dfrac{\partial \boldsymbol{\phi}}{\partial t}=\dot{\boldsymbol{\phi}}=\dfrac{\partial \boldsymbol{\phi}}{\partial t}(t)=\bigg\{\dfrac{\partial \phi_i}{\partial t}\bigg\}_{i=1}^{N_s}$, which are the vectors of control variables $\forall t \in (0, T)$ for the function $\phi_h$ and its time derivative $\dfrac{\partial \phi_h}{\partial t}$ respectively. Consider $[t_0 = 0, t_1, ..., t_{N_T} = T]$ of $N_T + 1$ equally spaced times in $[0, T]$. In this framework, we perform the time discretization by means of the generalized-$\alpha$ method \cite{Chung, Jensen, Tagliabue}: at time $t_j$, given $\dot{\boldsymbol{\phi}}_j$ and $\boldsymbol{\phi}_j$, find $\dot{\boldsymbol{\phi}}_{j+1}$, $\boldsymbol{\phi}_{j+1}$, $\dot{\boldsymbol{\phi}}_{j+\alpha_m}$, $\boldsymbol{\phi}_{j+\alpha_f}$ such that:
\begin{equation}
\begin{cases}
\boldsymbol{R}\left(\dot{\boldsymbol{\phi}}_{j+\alpha_m}, \boldsymbol{\phi}_{j+\alpha_f}\right) = \boldsymbol{0}, \\
\dot{\boldsymbol{\phi}}_{j+\alpha_m} = \dot{\boldsymbol{\phi}}_{j} + \alpha_m \left(\dot{\boldsymbol{\phi}}_{j+1} - \dot{\boldsymbol{\phi}}_{j}\right), \\
\boldsymbol{\phi}_{j+\alpha_f} = \boldsymbol{\phi}_{j} + \alpha_f (\boldsymbol{\phi}_{j+1}-\boldsymbol{\phi}_{j}), \\
\boldsymbol{\phi}_{j+1} = \boldsymbol{\phi}_{j} + \Delta t \dot{\boldsymbol{\phi}}_{j} + \delta \Delta t \left( \dot{\boldsymbol{\phi}}_{j+1} - \dot{\boldsymbol{\phi}}_{j} \right),
\label{eqn: NavierStokesStreamTimeDiscretization}
\end{cases}
\end{equation}
where $\delta, \alpha_m, \alpha_f \in \mathbb{R}_0^+$ are chosen on the basis of accuracy and stability considerations and $\Delta t$ is the fixed time step.

For further details about the numerical scheme and IGA we refer to \cite{Tagliabue}.

\section{Numerical results}
\label{section: NumericalResults}
We present some numerical results of wave equation and Navier-Stokes equations. For the first test case we compute all the variables at the mesh nodes, i.e. at the vertices, considering a mesh with $256$ elements. We use B–Splines basis functions of degree $p=2$ on uniform mesh of size $h = 1 / 64$ (4096 elements) and $\Delta t = 0.1$ for the lid-driven cavity problem. Indeed, IGA generally permits to obtain accurate numerical solutions using a reduced number of mesh elements than the Finite Element Method (FEM) \cite{Tagliabue}. This leads to smaller datasets that ease both the NN training and testing phases regardless of the specific ROM technique.

We implemented a MATLAB code for the numerical simulations of the wave equation, whereas the lid-driven cavity problem is available as a test case in the C++ IGA library isoglib.
We train the neural network using the Google Tesla K80 NVIDIA GPU.

We compute the KPOD-NN relative approximation error in the following way:
\begin{equation}
    \varepsilon_{KPOD-NN}(n, t, \boldsymbol{\mu}) = \dfrac{||\boldsymbol{u}_h(t; \boldsymbol{\mu}) - \mathbb{V} \boldsymbol{u}^{NN}_h(t; \boldsymbol{\mu}) ||}{|| \boldsymbol{u}_h(t; \boldsymbol{\mu}) ||}.
\end{equation}
This error is evaluated on a test dataset $\mathcal{D}_{test}$ made of $N_{test}$ elements. The final estimator that we use is the average of all KPOD-NN relative approximation errors computed on $N_{test}$ data:
\begin{equation}
    \overline{\varepsilon}_{KPOD-NN}(n) = \dfrac{\sum_{(t, \boldsymbol{\mu}) \in \mathcal{D}_{test}} \varepsilon_{KPOD-NN}(n, t, \boldsymbol{\mu})}{N_{test}}
\end{equation}

\subsection{Wave equation}

In this section, we focus on the numerical results related to the wave equation.
We denote the structure of both the training and test sets in Tables \ref{table: DatasetWave1} and \ref{table: DatasetWave2}. Values of $A_0$, $x_0$ and $\sigma$ in the building dataset are obtained by means of latin hypercube sampling \cite{Quarteroni2}.
We observe that, given a certain tolerance $\hat{\varepsilon}$ and a suitable value for $\gamma$, the reduced dimension that we get from KPOD ($n_{KPOD-NN}=15$) is fairly lower than the one obtained from POD ($n_{POD-NN}=109$). Moreover, the KPOD-NN size will be smaller than the POD-NN one and the former will be less expensive to train and to test in terms of both computational resources and total computational time.

\begin{table}[h!]
\centering
\begin{tabular}{|c c c c c c c c|}
    \hline
        Type of dataset & $N_h$ & $N_t$ & $N_{A_0}$ & $N_{x_0}$ & $N_{\sigma}$ & $n_{KPOD-NN}$ & $n_{POD-NN}$ \\ [0.5ex]
    \hline
        Training set & 256 & 100 & 5 & 5 & 5 & 15 & 109 \\
        Test set & 256 & 100 & 1 & 1 & 1 & - & - \\
    \hline
\end{tabular}
\vspace{10pt}
\caption{Building and test datasets for wave equation. These datasets have both the same mesh nodes and time steps, whereas values of $A_0$, $x_0$ and $\sigma$ are different. We consider $\gamma = 10^{-10}$ to determine $n_{KPOD-NN}$.}
\label{table: DatasetWave1}
\end{table}

\begin{table}[h!]
\centering
\begin{tabular}{|c c c|}
    \hline
        Type of dataset & Parameter & Values \\ [0.5ex]
    \hline
        \multirow{3}{*}{Training set} & $A_0$ & [0.5 0.641 0.721 0.821 1.0] \\
         & $x_0$ & [4.189 5.169 6.065 7.426 8.378] \\
         & $\sigma$ & [0.5 0.637 0.745 0.898 1.0] \\ \hline
        \multirow{3}{*}{Test set} & $A_0$ & 0.75 \\
         & $x_0$ & 8.0 \\
         & $\sigma$ & 0.9 \\
    \hline
\end{tabular}
\vspace{10pt}
\caption{Composition of building and test sets in terms of parameter values.}
\label{table: DatasetWave2}
\end{table}

\begin{figure}[t!]
	\centering
	\includegraphics[width=0.8\textwidth]{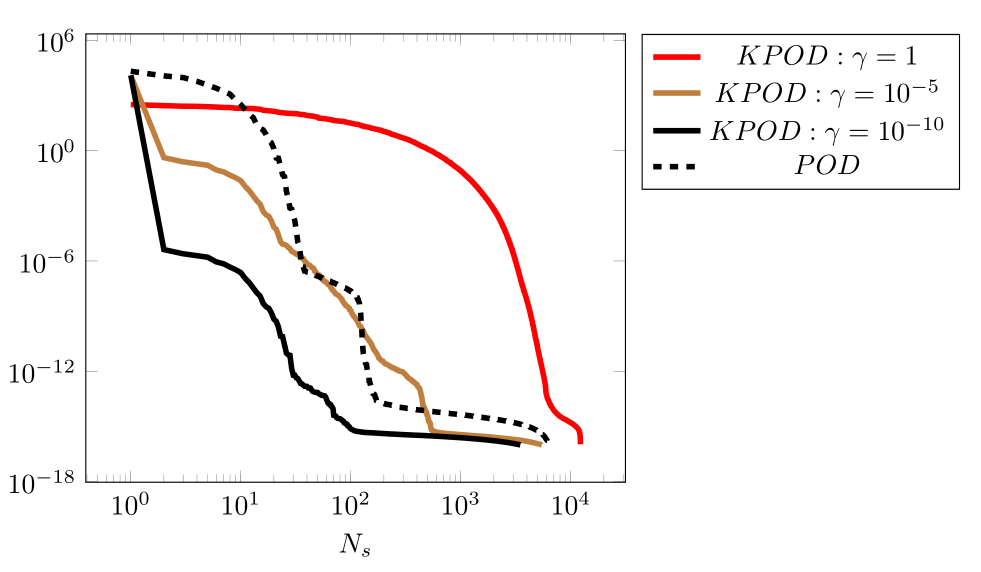}
	\caption{Comparison between the eigenvalues of the kernel matrix, coming from KPOD and different values for $\gamma$ (continuous), and the eigenvalues of the Gram matrix $\mathbb{C}$, coming from POD (dashed), for wave equation.}
	\label{fig: singularvaluesAlembert}
\end{figure}

\begin{figure}[t!]
	\centering
	\includegraphics[width=0.6\textwidth]{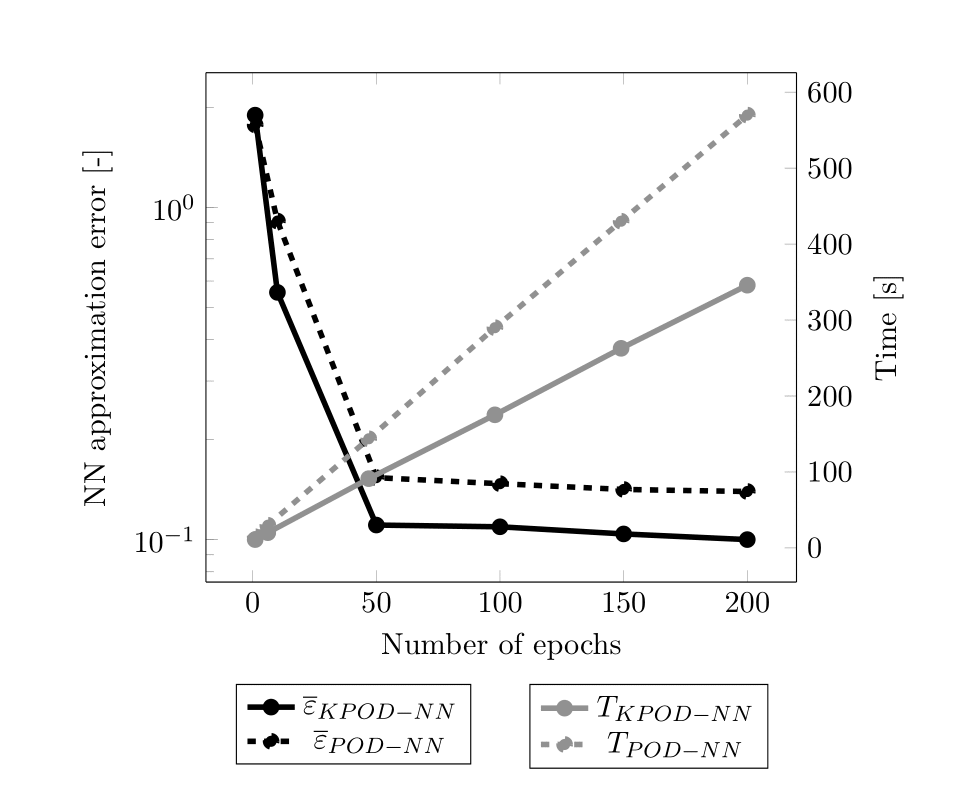}
	\caption{NN approximation errors and NN training computational times related to KPOD-NN and POD-NN vs. number of epochs for wave equation ($n_{KPOD-NN}$ = 15, $n_{POD-NN}$ = 109, $\gamma = 10^{-10}$).}
	\label{fig: NNApproximationErrorAlembert_1}
\end{figure}

\begin{figure}[t!]
	\centering
	\includegraphics[width=0.6\textwidth]{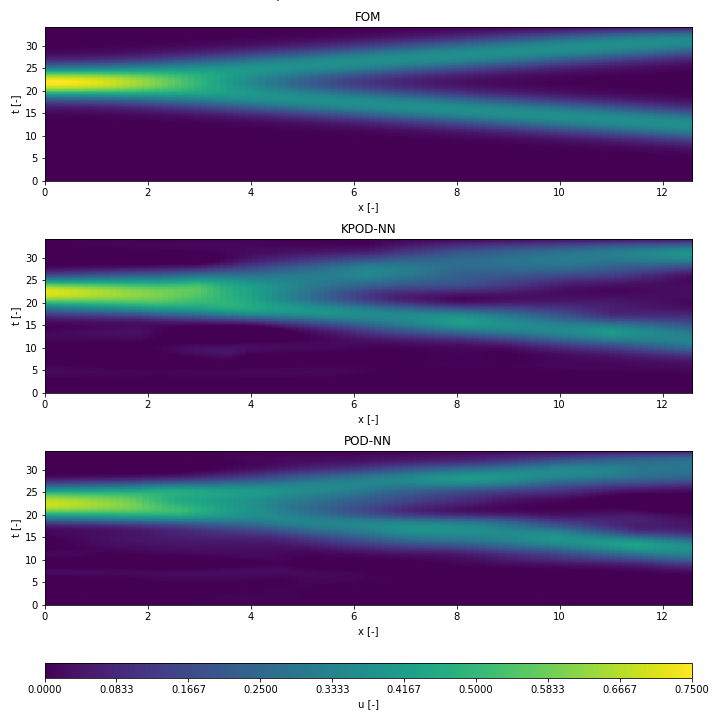}
	\caption{Numerical results on the test set $\mathcal{D}_{test}$ for wave equation. Comparison among FOM, KPOD-NN and POD-NN solutions. Both NN have been trained for $N_{epochs} = 100$, with $n_{KPOD-NN} = 15$ ($\gamma = 10^{-10}$) and $n_{POD-NN} = 109$.}
	\label{fig: numericalresultsAlembert}
\end{figure}

\begin{figure}[t!]
	\centering
	\includegraphics[width=0.6\textwidth]{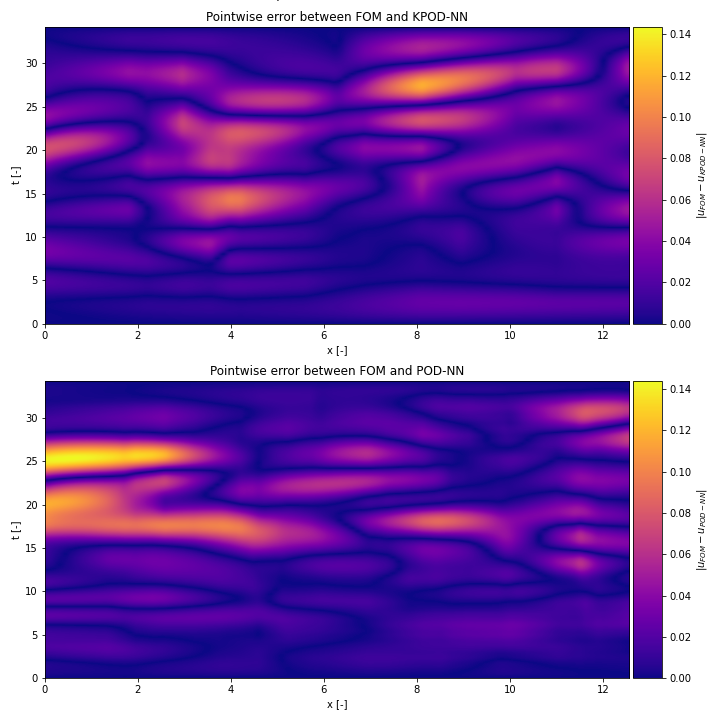}
	\caption{Pointwise absolute error on the test set $\mathcal{D}_{test}$ for wave equation. Comparison among FOM, KPOD-NN and POD-NN solutions. Both NN have been trained for $N_{epochs} = 100$, with $n_{KPOD-NN} = 15$ ($\gamma = 10^{-10}$) and $n_{POD-NN} = 109$.}
	\label{fig: errorsAlembert}
\end{figure}

\begin{figure}[t!]
	\centering
	\includegraphics[width=0.6\textwidth]{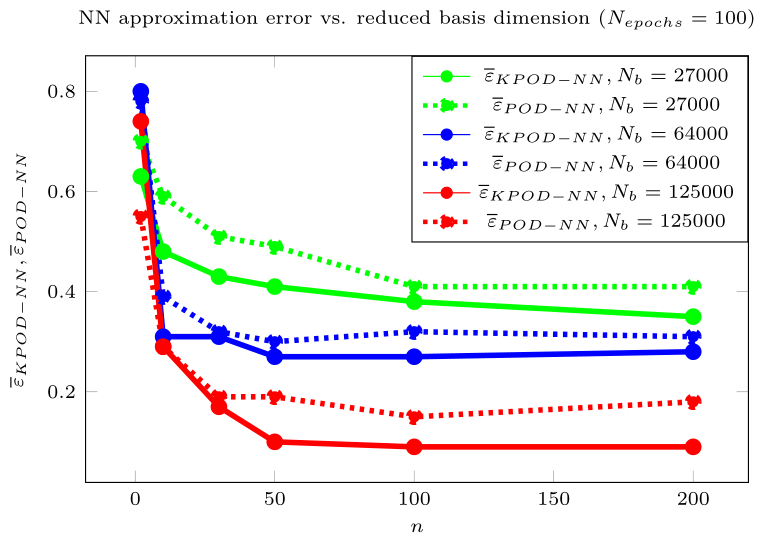}
	\caption{$\overline{\varepsilon}_{KPOD-NN}(n)$ and $\overline{\varepsilon}_{POD-NN}(n)$ errors vs. $n$ for different choices of the dataset, which contain the numerical solutions of the wave equation for different parameters instances. $N_b = 125000$ corresponds to the building set introduced in Tables \ref{table: DatasetWave1} and \ref{table: DatasetWave2}.}
	\label{fig: NNApproximationErrorAlembert_2}
\end{figure}

In Figure \ref{fig: singularvaluesAlembert} we highlight the role of $\gamma$. In particular, at least for this test case and for the range of values that we consider, the lower the value of parameter $\gamma$ is, the higher the rate of decay is. We show the capability of the KPOD method to potentially collect the most significant modes and information in the really first eigenvalues of the kernel matrix.
Moreover, the computational time that is needed to perform KPOD is comparable with the one of POD, due to the fact that we do not have to deal with high-dimensional data explicitly \cite{Scholkopf}.

We depict in Figure \ref{fig: NNApproximationErrorAlembert_1} the behavior of both $\overline{\varepsilon}_{KPOD-NN}(n)$ and $\overline{\varepsilon}_{POD-NN}(n)$ errors with respect to $N_{epochs}$. We see that, providing a proper number of epochs, once convergence is reached, the approximation error related to KPOD-NN is slightly smaller than the one of POD-NN method, even if the size of the two NN is different. We also notice that the total time to train the KPOD-NN is approximately halved with respect the POD-NN one.
By looking at the numerical results in Figure \ref{fig: numericalresultsAlembert} on the test set $D_{test}$ for all spatial and temporal coordinates, we see a good agreement between the FOM solution and the KPOD-NN one.
This is also confirmed in Figure \ref{fig: errorsAlembert}, where we show that also the pointwise difference in absolute value between the FOM solution and the ROM one is generally lower for the KPOD-NN method.
In Figure \ref{fig: NNApproximationErrorAlembert_2} we see the plot related to NN approximation error vs. reduced dimension $n$ for different building sets. We comment that, given a specific NN architecture and size, KPOD-NN has a higher representational power than POD-NN, and the smaller the training dataset is, the higher the plateau of the approximation error is.
We stress that in this specific test case we are solving a linear PDE in 1D.
For this reason, the advantages of KPOD over POD on the NN approximation errors are quite limited.

\begin{figure}[t!]
	\centering
	\includegraphics[width=1.0\textwidth]{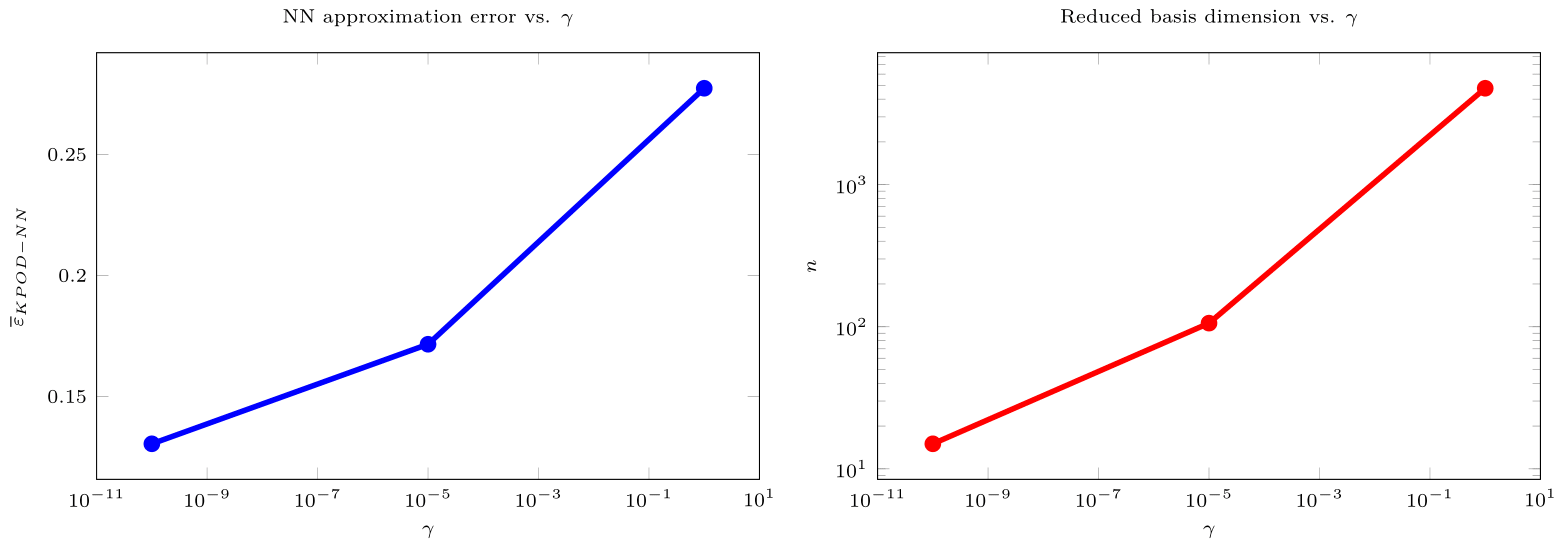}
	\caption{NN approximation error $\overline{\varepsilon}_{KPOD-NN}$ (left) and reduced basis dimension $n$ (right) vs. $\gamma$ for wave equation. All NN have been trained for $N_{epochs}$ = 100.}
	\label{fig: epsilon_KPOD_NN}
\end{figure}

\begin{figure}[t!]
	\centering
	\includegraphics[width=1.0\textwidth]{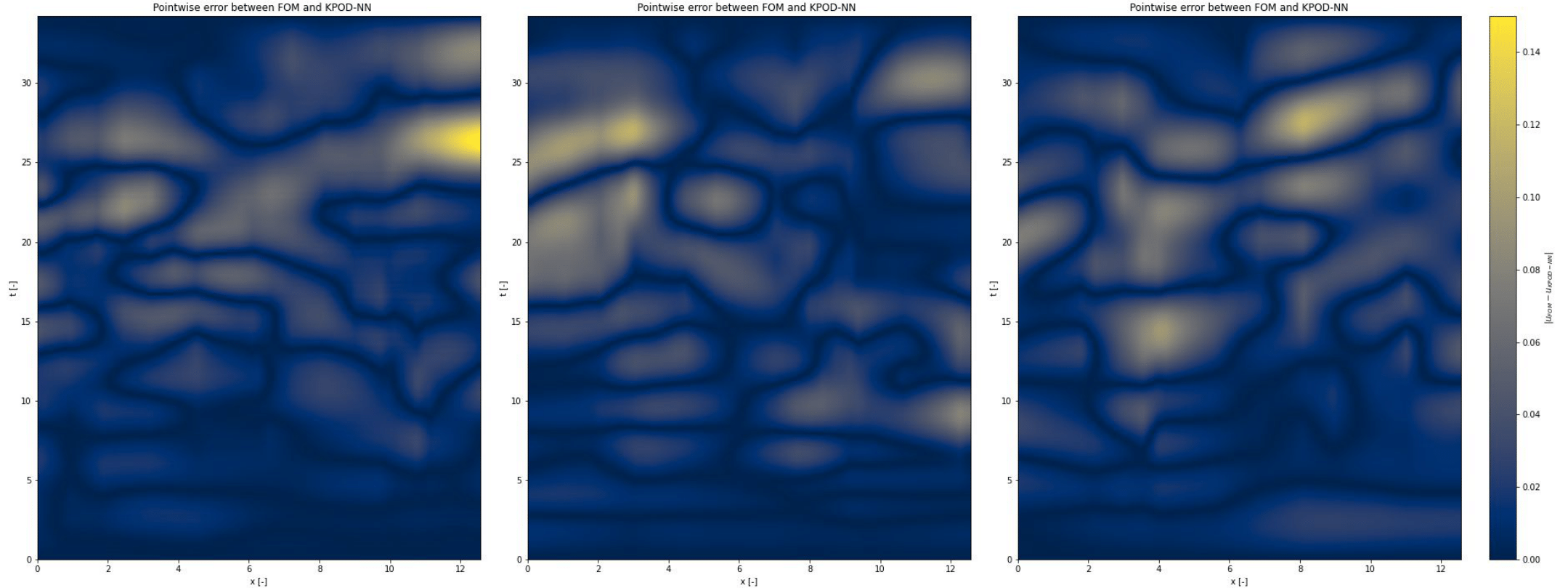}
	\caption{Pointwise absolute error on the test set $\mathcal{D}_{test}$ for wave equation. Comparison among KPOD-NN solutions for different values of $\gamma$. All NN have been trained for $N_{epochs} = 100$. We use $\gamma=1$ (left, $n_{KPOD-NN} = 4772$), $\gamma=10^{-5}$ (center, $n_{KPOD-NN} = 106$) and $\gamma=10^{-10}$ (right, $n_{KPOD-NN} = 15$).}
	\label{fig: test_Alembert_error_gamma}
\end{figure}

In Figure \ref{fig: epsilon_KPOD_NN} we show both the behavior of the NN approximation error $\overline{\varepsilon}_{KPOD-NN}$ and reduced basis dimension $n$ with respect to parameter $\gamma$. Given a fixed tolerance $\hat{\varepsilon}$, we notice that both $\overline{\varepsilon}_{KPOD-NN}$ and $n$ scales monotonically with $\gamma$. We also depict in Figure \ref{fig: test_Alembert_error_gamma} the pointwise absolute error related to KPOD-NN solutions for different choices of $\gamma$. We highlight that high values, such as $\gamma=1$, are not favorable and might be source of high localized errors in space and time. Indeed in this scenario, given again a fixed tolerance $\hat{\varepsilon}$, the reduced dimension $n$ could be potentially high and could lead to a bigger NN, which is more challenging to train. On the other hand, when lower values are considered (e.g. $\gamma=10^{-5}$ or $\gamma=10^{-10}$), KPOD modes decay faster and steeply while the corresponding NN is smaller and easier to train.

\subsection{Lid-driven cavity benchmark}

\begin{figure}[t!]
	\centering
	\includegraphics[width=0.6\textwidth]{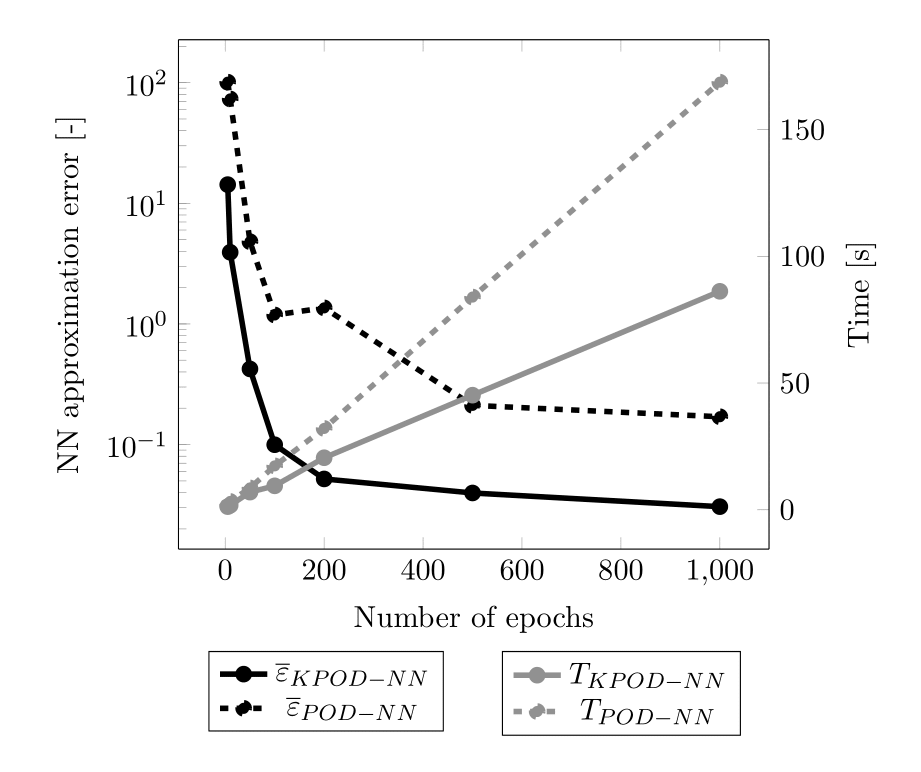}
	\caption{NN approximation errors and NN training computational times related to KPOD-NN and POD-NN vs. number of epochs for the lid driven cavity benchmark ($n_{KPOD-NN}$ = 12, $n_{POD-NN}$ = 131).}
	\label{fig: NNApproximationErrorLidDrivenCavity}
\end{figure}

Here we present the numerical results for the lid-driven cavity problem.
We generate the training set by means of FOM solutions, considering evenly spaced Reynolds numbers between 100 and 5000 with a fixed step equal to 10.
The test set is sampled for $\mathbb{R}e \in [105, 4995]$ by using again a fixed step equal to 10.
In Table \ref{table: PerformancesLiddrivencavity} and Figure \ref{fig: NNApproximationErrorLidDrivenCavity} we report some information related to the training phase of the NN when either KPOD or POD is employed on the same building set. Given a certain number of epochs and a tolerance $\hat{\varepsilon}$, the KPOD-NN approach leads to an approximation error that is one order of magnitude less than the POD-NN one. Moreover, the computational time is again halved for KPOD-NN and its reduced basis dimension is reduced by a factor 10. We highlight that the dimension of the reduced basis affects directly the number of NN parameters to train, i.e. the weights of the NN. This number is significantly smaller for KPOD-NN.

\begin{table}[h!]
\centering
\begin{tabular}{|c c c c c c|}
    \hline
        Technique & NN parameters & $n$ & Number of epochs & NN approximation error & Computational time $[s]$ \\ [0.5ex]
    \hline
        KPOD-NN & 360 & 12 & 1000 & $3.1 \cdot 10^{-2}$ & 86.17 \\
        POD-NN & 69'954 & 131 & 1000 & $1.7 \cdot 10^{-1}$ & 168.7 \\
    \hline
\end{tabular}
\vspace{10pt}
\caption{Comparison of the KPOD-NN and POD-NN methods for the lid-driven cavity problem on the same building set, in terms of both $\overline{\varepsilon}_{KPOD-NN}(n)$/$\overline{\varepsilon}_{POD-NN}(n)$ (i.e. NN approximation error) and computational times to train the NN. We use $\gamma=10^{-5}$ for KPOD.}
\label{table: PerformancesLiddrivencavity}
\end{table}

\begin{figure}[h!]
	\centering
	\includegraphics[width=0.5\textwidth]{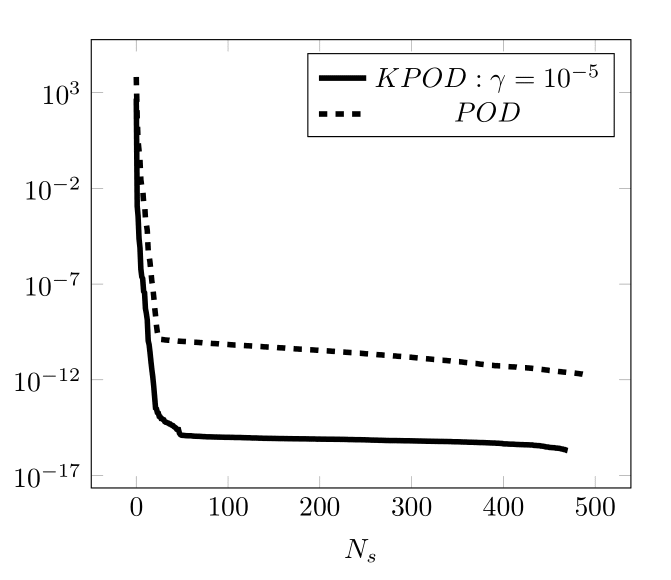}
	\caption{Comparison between the eigenvalues of the kernel matrix, coming from KPOD (continuous), and the eigenvalues of the Gram matrix $\mathbb{C}$, coming from POD (dashed), for the lid-driven cavity problem.}
	\label{fig: singularvaluesliddrivencavity}
\end{figure}

Moreover, also in this case, we get a reduced value of $n$ thanks to the better behavior of the eigenvalues of the kernel matrix with respect to the singular values of matrix $\mathbb{S}$, as can be seen in Figure \ref{fig: singularvaluesliddrivencavity}. The computational times to train the NN coming from the KPOD approximation is strongly reduced with respect to the one related to POD. Testing times remain comparable and negligible. The same occurs for the computational times to perform KPOD and POD effectively.

\begin{figure}[t!]
\centering
\begin{subfigure}{0.40\textwidth}
  \centering
  \includegraphics[width=\textwidth]{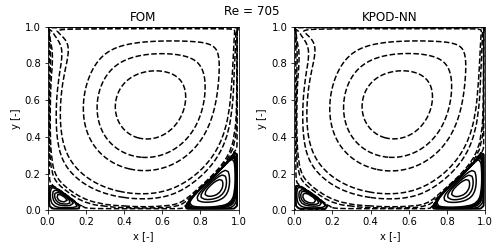}
\end{subfigure}%
~
\begin{subfigure}{0.40\textwidth}
  \centering
  \includegraphics[width=\textwidth]{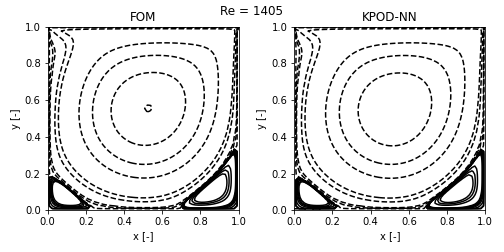}
\end{subfigure}

\begin{subfigure}{0.40\textwidth}
  \centering
  \includegraphics[width=\textwidth]{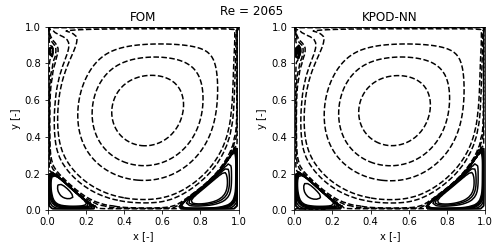}
\end{subfigure}%
~
\begin{subfigure}{0.40\textwidth}
  \centering
  \includegraphics[width=\textwidth]{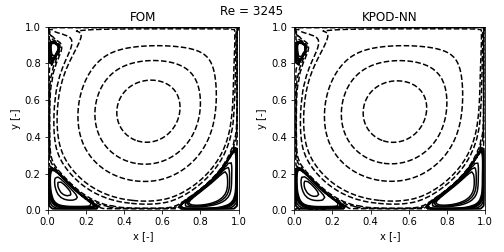}
\end{subfigure}

\begin{subfigure}{0.40\textwidth}
  \centering
  \includegraphics[width=\textwidth]{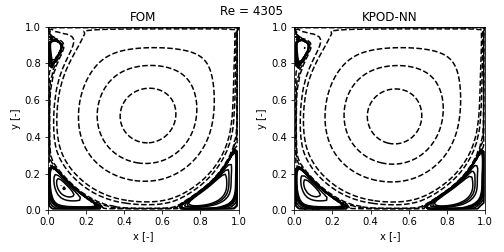}
\end{subfigure}%
~
\begin{subfigure}{0.40\textwidth}
  \centering
  \includegraphics[width=\textwidth]{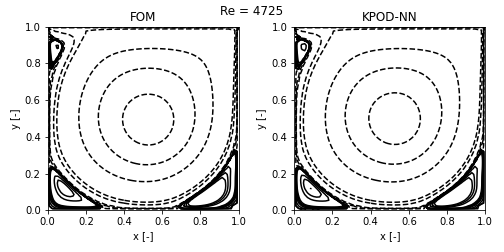}
\end{subfigure}

\caption{Streamlines computed on the test set $\mathcal{D}_{test}$ for the lid-driven cavity benchmark. Comparison between FOM and KPOD-NN solutions. The NN has been trained for $N_{epochs} = 1000$, with $n_{KPOD-NN} = 12$ ($\overline{\varepsilon}_{KPOD-NN} = 3.1 \cdot 10^{-2}$).}
\label{fig: numericalresultsliddrivencavitystreamlines}
\end{figure}

\begin{figure}[h!]
\centering
\begin{subfigure}{0.45\textwidth}
  \centering
  \includegraphics[width=\textwidth]{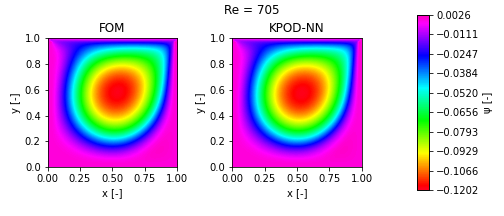}
\end{subfigure}%
~
\begin{subfigure}{0.45\textwidth}
  \centering
  \includegraphics[width=\textwidth]{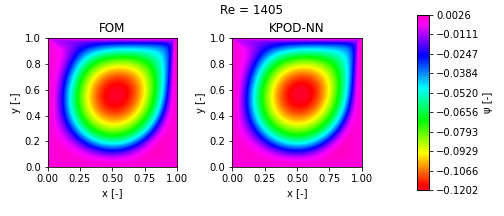}
\end{subfigure}

\begin{subfigure}{0.45\textwidth}
  \centering
  \includegraphics[width=\textwidth]{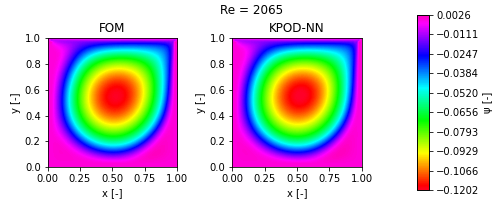}
\end{subfigure}%
~
\begin{subfigure}{0.45\textwidth}
  \centering
  \includegraphics[width=\textwidth]{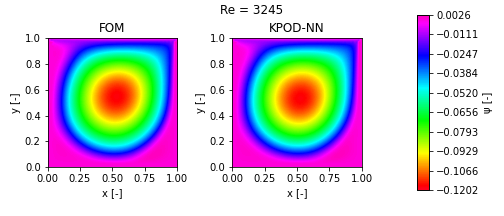}
\end{subfigure}

\begin{subfigure}{0.45\textwidth}
  \centering
  \includegraphics[width=\textwidth]{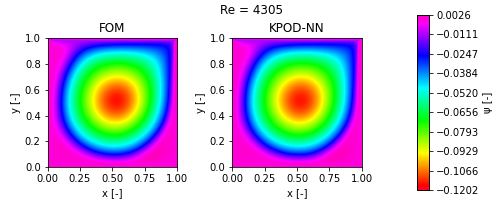}
\end{subfigure}%
~
\begin{subfigure}{0.45\textwidth}
  \centering
  \includegraphics[width=\textwidth]{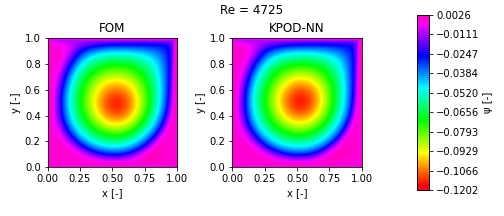}
\end{subfigure}

\caption{Streamfunction $\phi$ computed on the test set $\mathcal{D}_{test}$ for the lid-driven cavity benchmark. Comparison between FOM and KPOD-NN solutions. The NN has been trained for $N_{epochs} = 1000$, with $n_{KPOD-NN} = 12$ ($\overline{\varepsilon}_{KPOD-NN} = 3.1 \cdot 10^{-2}$).}
\label{fig: numericalresultsliddrivencavitystreamfunction}
\end{figure}

In Figures \ref{fig: numericalresultsliddrivencavitystreamlines} and \ref{fig: numericalresultsliddrivencavitystreamfunction} we compare the streamlines and the values assumed by the streamfunction $\phi$ for FOM and KPOD-NN for different Reynolds numbers $\mathbb{R}e$. We see that, considering the number and positions of the contour lines, the KPOD-NN method captures all vortexes properly. Indeed strong nonlinearities arises in the 2D lid-driven cavity problem, where there are both a discontinuous dependance of the solution with the Reynolds number and bifurcations in the development of vortices. For these reasons, in this second test case, we have stronger advantages in using a fully nonlinear technique, such as KPOD-NN.

~

\section{Conclusions}
\label{section: Conclusions}
In this work, we proposed a non-intrusive reduced order model technique that combines KPOD with an adaptively built NN, whose number of layers and number of neurons scale according to the dimension of the extracted reduced basis. The use of nonlinear dimensionality reduction in the first part of the algorithm permits to shrink the valuable information on the first modes. This operation strongly reduces the number of functions that must be collected to obtain a basis up to specific tolerance, in particular if this method is compared to linear dimensionality reduction techniques, such as POD.

We applied our methodology to parametrized parabolic and hyperbolic PDEs in both 1D and 2D settings, either in a linear (wave equation) or nonlinear (Navier-Stokes equations) context. We saw a good agreement between the FOM solution and the KPOD-NN solution in both cases. We also compared our method with the POD-NN one, by showing that the reduced coefficients coming from KPOD contain potentially more information than the one extracted by means of POD. This leads to a smaller NN approximation error on the test sets. Moreover, we remind that the NN size coming from KPOD is again by construction necessarily smaller and easier to train than the one given by POD. Both KPOD and POD involve the same computational resources, and the computational times spent by the two methods to compute the reduced basis is still similar. According to our tests, by fixing a priori a certain tolerance, KPOD-NN collects a number of modes which is at least 10 times smaller than POD-NN while leading to smaller NN approximation errors, which are reduced by approximately a factor of 1.2 for the wave equation and by approximately a factor of 10 for the lid-driven cavity problem, where nonlinearities arise. Moreover, the training costs for KPOD-NN are at least halved with respect to POD-NN ones. By using even lower values of $\gamma$, all the advantages of KPOD-NN over POD-NN are potentially even stronger.

A possible extension of this work could address 3D problems with a significantly higher number of DOFs. Another possible topic for this non-intrusive technique could be the application to multifield fully-coupled problems, such as the one arising in cardiac modeling: in this framework, the nonlinear PDEs related to electrophysiology, mechanics and fluid dynamics can be reduced in an independent manner with potentially different NNs. Finally, we notice that our KPOD-NN can also be embedded in pre-existing deep learning ROM frameworks to improve their efficiency and accuracy.

\section*{Acknowledgements}
We sincerely thank Prof. A. Quarteroni, Dr. F. Regazzoni and Dr. S. Fresca for the useful discussions about numerical analysis, nonlinear dimensionality reduction and deep learning.

\printbibliography

\end{document}